\documentclass[11pt]{amsart}
\usepackage{amsmath,amsfonts}
\usepackage{amssymb}

\usepackage[all]{xypic}

\usepackage{amscd}

\newtheorem{thm}{Theorem}[subsection]
\newtheorem{lem}[thm]{Lemma}
\newtheorem{prop}[thm]{Proposition}

\newtheorem{cor}[thm]{Corollary}
\newtheorem{rmk}[thm]{Remark}

\newtheorem{ex}[thm]{Example}

\def\del{\partial}

\def\R{\mathbb{R}}
\def\C{\mathbb{C}}
\def\Z{\mathbb{Z}}
\def\P{\mathbb{P}}

\def\on{\operatorname}
\def\cD{\mathcal D}
\def\CD{\mathcal D}
\def\CF{\mathcal F}

\def\CK{\mathcal K}

\def\CN{\mathcal N}
\def\CO{\mathcal O}

\def\CS{\mathcal S}
\def\CT{\mathcal T}
\def\Hom{{\on{Hom}}}
\def\ol{\overline}
\def\intHom{{\mathcal Hom}}

\def\CY{\mathcal Y}

\def\hom{{hom}}

\def\CC{\mathcal C}

\def\fU{\mathfrak U}
\def\CM{\mathcal M}

\def\CE{\mathcal E}

\def\fY{\mathfrak Y}
\def\fS{\mathfrak S}

\def\codim{\on{codim}}

\def\ul{\underline}
\def\id{{\it{id}}}

\def\Vect{{\it{Vect}}}

\def\Top{{\it{Top}}}
\def\Ch{{\it{Ch}}}
\def\Sh{{\it{Sh}}}

\def\hra{\hookrightarrow}

\def\orient{{\it{or}}}

\def\Ext{{\on{Ext}}}

\def\ff{\mathfrak f}

\def\risom{\stackrel{\sim}{\to}}
\def\rmod{mod_r}
\def\lmod{mod_\ell}

\def\CP{\mathcal P}

\def\fra{\mathfrak a}
\def\frb{\mathfrak b}

\def\frI{\mathfrak I}

\def\fraki{\mathfrak i}

\def\fL{\mathfrak L}

\begin{document}

\title[Microlocal branes are constructible sheaves]
{Microlocal branes are constructible sheaves}

\author[David Nadler]{David Nadler\\
\\
\it{\lowercase{\uppercase{D}edicated to \uppercase{P}aul \uppercase{S}. \uppercase{N}adler}}}

\begin{abstract} 
Let $X$ be a compact real analytic manifold, and let $T^*X$ be its cotangent bundle.
In a recent paper with E. Zaslow~\cite{NZ}, we showed that the dg category $Sh_c(X)$ of
constructible sheaves on $X$ quasi-embeds
into the triangulated envelope $F(T^*X)$
of the Fukaya category of $T^*X$.
We prove here that the quasi-embedding is in fact a quasi-equivalence.
When $X$ is a complex manifold,
one may interpret this as a topological analogue of the identification of
Lagrangian branes in $T^*X$ and regular holonomic $\cD_X$-modules developed
by Kapustin~\cite{Kap} and
Kapustin-Witten~\cite{KW} from a physical perspective.

As a concrete application, we show that compact connected
exact Lagrangians in $T^*X$ (with some modest homological assumptions) 
are equivalent in the Fukaya category to the zero section.
In particular, this determines their (complex) cohomology ring and homology
class in $T^*X$, and provides a homological bound on their number of intersection points. 
An independent characterization of compact branes in $T^*X$ has recently been 
obtained by Fukaya-Seidel-Smith~\cite{FSS}.
\end{abstract}

\maketitle


{\tiny
\tableofcontents
}


\section{Introduction}

\subsection{Summary}

Let $X$ be a compact real analytic manifold, and let $T^*X$ be its cotangent bundle.
Let $Sh_c(X)$ be the differential graded (dg) category of constructible complexes of sheaves
on $X$. Objects of $Sh_c(X)$ are complexes of sheaves of $\C$-vector spaces on $X$
with bounded constructible cohomology; morphisms are obtained from the naive 
morphism complexes by inverting quasi-isomorphisms.
The (ungraded) cohomology category of $Sh_c(X)$ is the usual bounded derived category $D_c(X)$
of cohomologically constructible complexes on $X$. (See Section~\ref{sec sheaves} below
for more details.)

In a joint paper with Eric Zaslow~\cite{NZ} (reviewed in Section~\ref{sec fukaya} below),
we developed a Fukaya $A_\infty$-category of exact Lagrangian branes in $T^*X$.
Objects are (not necessarily compact) exact Lagrangian submanifolds of $T^*X$
equipped with brane structures and perturbations; 
morphisms are given by transverse intersections,
and composition maps by counts of pseudoholomorphic polygons.
A key assumption on the objects is that they have reasonable compactifications
so that we can make sense of ``intersections at infinity".
Perturbations are organized so that morphisms always propagate 
a small amount ``forward in time". 
The construction is a close relative of the category of vanishing cycles 
proposed by Kontsevich~\cite{Kontlecture} and Hori-Iqbal-Vafa~\cite{HIV},
and developed by Seidel~\cite{SeidmutationI}, \cite{SeidmutationII}, \cite{Seidel}.

This version of the Fukaya category of the cotangent bundle
encodes the geometry of infinitesimal paths in the base.
Many other theories (such as 
cyclic homology~\cite{Loday}, and 
the chiral de Rham complex~\cite{KV})
also relate small loops in the base to the de Rham complex (and hence $\CD$-modules),
while other versions of Floer theory on the cotangent bundle 
are closely related to the full path space of the base (for example, in the work
of Viterbo~\cite{VitICM}).

\medskip

Now let us pass to a stable setting and consider the dg category of right modules over the 
above Fukaya category. 
We write $F(T^*X)$ for
the full subcategory
of twisted complexes of representable modules,
and refer to it  as the triangulated envelope of the Fukaya category.
We write $DF(T^*X)$ for the cohomology category of $F(T^*X)$, and 
refer to it as the derived Fukaya category.

The main result of~\cite{NZ} was the construction of 
an $A_\infty$-quasi-embedding
$$
\xymatrix{
\mu_X:Sh_c(X) \ar@{^(->}[r] & F(T^*X)
}
$$
which we will refer to as microlocalization.
This is an $A_\infty$-functor such that when we pass to cohomology, we obtain a fully faithful 
embedding of 
triangulated categories
$$
\xymatrix{
H(\mu_X):D_c(X)\ar@{^(->}[r] & DF(T^*X).
}
$$

Our aim in this paper is to show that $\mu_X$ is in fact a quasi-equivalence. This is the assertion
that the embedding $H(\mu_X)$ of triangulated categories
is an equivalence. 
In other words, we must see that it is essentially surjective, or in plain words that
$$
\mbox{\em{every object of $DF(T^*X)$
is isomorphic to an object coming from $D_c(X)$.}}
$$

This can be viewed as a categorification of the fact~ \cite[Theorem 9.7.10]{KS} that the characteristic
cycle homomorphism is an isomorphism from constructible functions to conical Lagrangian cycles.
To simplify the statement, let us assume for the moment that $X$ is oriented.
Then as a consequence of our main result, one can deduce a commutative diagram
relating the microlocalization $\mu_X$ to the characteristic cycle homomorphism
$$
\xymatrix{
K_0(D_c(X)) \ar[rr]^-{K_0(H(\mu_X))} \ar[d]_-{\chi}^-{\wr} & & \ar[d]^-{\xi} K_0(DF(T^*X))\\
\mathfrak{F}_c(X) \ar[rr]^-{CC} & & \mathfrak L_{con}({T^*X})
}
$$
Here $K_0$ denotes the underlying Grothendieck group,
$\mathfrak{F}_c(X)$ the group of constructible functions on $X$, and
$\mathfrak{L}_{con}({T^*X})$ the group of conical Lagrangian cycles.
The homomorphism $\chi$ is the
local Euler-Poincar\'e index, and the homomorphism $\xi$ simply dilates a Lagrangian brane
down to a conical cycle. 
The reader could consult~\cite{KS} for a comprehensive treatment of the 
characteristic cycle homomorphism $CC$, local Euler-Poincar\'e index $\chi$
and related topics.
The fact that the diagram is commutative follows
from the compatibility of the definition of $\mu_X$ with the functoriality formula~\cite[Theorem 4.2]{SV}
(see also the results of~\cite{Gins} in the complex algebraic setting).


\subsection{Sketch of arguments}

To establish our main assertion, we use the following variation of a standard argument.
Consider the toy problem of showing a set of vectors $\{v_\alpha\}$ span a 
finite dimensional vector space $V$.
We can reformulate this in terms of the identity map $\id_V:V\to V$.
Namely, we can ask if $\id_V$ can be written in the form
$$
\id_V=\sum_\alpha \lambda_\alpha\otimes v_\alpha
$$ 
where $\{\lambda_\alpha\}$ are some set of functionals. If so, then 
by applying $\id_V$ to any $v\in V$, we see that $v$ is in the span of 
the set $\{v_\alpha\}$.
We may interpret this argument as expressing the identity map $id_V$
as a sum of projections $ \lambda_\alpha\otimes v_\alpha$ onto the span of our collection.
More generally, if $V$ is not necessarily finite dimensional, it still suffices to show
that for any finite dimensional subspace $V'\subset V$, the restriction of the identity
$\id_V|_{V'}$ can be expressed in the above form.

The same proof can
be implemented in the less elementary setting of a triangulated $A_\infty$-category $\CC$. Suppose we want to show that a collection of objects $\{c_\alpha\}$
classically generates all of $\CC$. In other words, we want to see that any object $c$ can be 
realized as a finite sequence of iterated cones beginning with 
maps among the objects $\{c_\alpha\}$.
Consider the identity functor
$\id_\CC$ in the 
triangulated $A_\infty$-category of $A_\infty$-functors from $\CC$ to $\CC$.
Then it suffices to show that $\id_\CC$ itself can be 
realized as a finite sequence of iterated cones of 
functors of the form $f_\alpha\otimes c_\alpha$, where $f_\alpha$ is a bounded object
of the dg category 
of left $A_\infty$-modules. 
In some situations such as the one we will consider below, 
the above is stronger than what can actually be shown. But
it still suffices to verify the weaker assertion
that for any full subcategory $\CC'$ of $\CC$ generated by finitely many objects,
the restriction $\id_\CC|_{\CC'}$
is in the full subcategory generated by functors of the form $ f_\alpha\otimes c_\alpha$.

As an example, consider the problem of representing coherent sheaves
on projective space $\P^n$ by complexes of vector bundles. 
In~\cite{Beil}, Beilinson introduced
a resolution of the structure sheaf $\CO_{\Delta_{\P^n}}$ of the diagonal $\Delta_{\P^n}
\subset \P^n\times \P^n$
by external products of vector bundles. For any coherent sheaf on $\P^n$,
convolution with this resolution produces the desired complex of vector bundles.

\medskip

We will apply the above strategy in the context of $F(T^*X)$ to see that every object comes from $Sh_c(X)$ via the microlocalization functor $\mu_X$. 

To put this plan
in action, we must have a way to get our hands
on functors between categories of branes, in particular the identity functor.
For this, we take advantage of the symmetry of cotangent bundles
and introduce a duality equivalence
$$
\alpha_X: F(T^*X)^\circ \to F(T^*X).
$$
By definition, it acts on the underlying Lagrangians of our branes by the antipodal
anti-symplectomorphism
$$
a:T^*X\to T^*X 
\qquad
a(x,\xi) = (a,-\xi).
$$
As a consequence of our main result, we will see that the brane duality $\alpha_X$
corresponds to Verdier duality $\CD_X$
under the microlocalization $\mu_X$.

Now given compact real analytic manifolds $X_0,X_1$,
we will construct functors 
$$
F(T^*X_0)\to F(T^*X_1)
$$ 
by thinking of branes in 
the product $T^*X_0\times T^* X_1$ as integral kernels.
Consider the dg category $\rmod(F(T^*X_1))$ of right $A_\infty$-modules
and the corresponding Yoneda embedding
$$
\CY_r:F(T^*X_1)\to \rmod(F(T^*X_1))
$$
$$
\CY_r(P_1): P'_1 \mapsto \hom_{F(T^*X_1)}(P'_1, P_1).
$$
For each object $L$ of $F(T^*X_0\times T^*X_1)$,
we define an $A_\infty$-functor by considering the mapping functional
$$
\tilde\Psi_{L*}: F(T^*X_0) \to  \rmod(F(T^*X_1))
$$
$$
\tilde\Psi_{L*}(P_0) :P_1\mapsto  \hom_{F(T^*X_0\times T^*X_1)}(L, P_0 \times \alpha_{X_1}(P_1)).
$$
Note that $\tilde\Psi_{L*}$ is functorial in $L$ in the contravariant sense.
As a consequence of our main result, we will see that there is a functor
$$
\Psi_{L*}: F(T^*X_0) \to  F(T^*X_1)
$$
that represents $\tilde\Psi_{L*}$ in the sense that we have a quasi-isomorphism of functors
$$
\tilde\Psi_{L*}\simeq \CY_r\circ\Psi_{L*} : F(T^*X_0) \to  \rmod(F(T^*X_1)).
$$

Two special cases of the above construction merit special mention.
First, for the microlocalization $L_\Delta =\mu_X(\C_{\Delta_X})$
of the constant sheaf $\C_{\Delta_X}$ along the diagonal $\Delta_X\subset X\times X$,
we check directly that the functor $\tilde \Psi_{L_\Delta *}$ 
is represented by the identity functor $\id_{F(T^*X)}$
in the sense that there is a quasi-isomorphism of functors
$$
\tilde\Psi_{L_\Delta *}\simeq \CY_r  : F(T^*X) \to  \rmod(F(T^*X)).
$$
Second, for an external product $L_0\times L_1$,
we check directly that the functor $\tilde\Psi_{L_0\times L_1}$
plays the expected role of a projection
in the sense that there is a quasi-isomorphism of functors
$$
\tilde\Psi_{L_0\times L_1*} \simeq 
\CY_\ell(L_0)
\otimes \CY_r(\alpha_{X_1}(L_1))
 : F(T^*X_0) \to  \rmod(F(T^*X_1)).
$$
Here 
$
\CY_\ell 
$
denotes the Yoneda embedding 
$$
\CY_\ell:F(T^*X_0)\to \lmod(F(T^*X_0))
$$
$$
\CY_\ell(P_0): P'_0 \mapsto \hom_{F(T^*X_0)}(P_0, P'_0)
$$
for left $A_\infty$-modules over $F(T^*X_0)$.
Note that the above formula says that $\tilde\Psi_{L_0\times L_1*} $ is in fact a projection onto
the span of the dual brane $\alpha_{X_1}(L_1)$. Although one might prefer
simpler formulas, our conventions are set up
to agree with  the pushforward of sheaves.

An alternative, more geometric 
framework for constructing functors between categories of branes is provided by 
the beautiful formalism of world-sheet foam introduced by Khovanov-Rozansky
from a physical viewpoint~\cite{KR},
or that of quilted Riemann surfaces developed
by Wehrheim-Woodward in a mathematical context~\cite{WW}. 
The latter was an original inspiration for 
the strategy of proof undertaken here.
To keep this paper as self-contained as possible, we have
opted for the above homological approach, though we have included a brief discussion explaining
its compatibility with that of~\cite{WW}. 

\medskip

Now to see that $F(T^*X)$ is classically generated by objects coming from $Sh_c(X)$,
we would like to realize the identity functor $id_{F(T^*X)}$
as an iterated cone of
external products $ f_\alpha \otimes \mu_X(\CF_\alpha)$, for some objects $\CF_\alpha$
of $Sh_c(X)$,
and some bounded objects $f_\alpha$ of $\lmod(F(T^*X_0))$.
It is not difficult to see that this strong an assertion can not be true.
Instead, we fix a conical Lagrangian $\Lambda\subset T^*X$ and consider
the full $A_\infty$-subcategory 
$$
F(T^*X)_\Lambda
\subset  
F(T^*X)
$$
consisting
of branes whose boundaries at infinity lie in the boundary of $\Lambda$.
By~\cite{NZ}, every finite collection of objects of $F(T^*X)$
lies in such a subcategory $F(T^*X)_\Lambda$ for some $\Lambda$.
Thus to arrive at our desired conclusion,
it suffices to realize the identity functor of $F(T^*X)_\Lambda$
as an iterated cone of
external products $f_\alpha \otimes \mu_X(\CF_\alpha)$. 

Translating this into the above setting of branes in the product $T^*X\times T^*X$,
we seek to express the restriction of the functor
$
\tilde \Psi_{L_\Delta *} 
$
to the subcategory $F(T^*X)_\Lambda$
as an iterated cone of the restrictions of functors of the form
$
\tilde \Psi_{L_\alpha \times \alpha_X(\mu_X(\CF_\alpha))}. 
$
By the functoriality of our constructions,
this would follow immediately if the brane $L_{\Delta_X}$
could be realized as an iterated cone
of branes of the form $L_\alpha \times \alpha_X(\mu_X(\CF_\alpha))$.
Of course, this is not true (for example, it would imply the identity functor $id_{F(T^*X)}$
in fact could be written 
in terms of external products), but we can achieve the following:
there is a 
collection of $[0,1]$-families of objects $\fL_{\Delta_{X_\alpha}, t}$ 
of $F(T^*X\times T^*X)$ satisfying the following properties:

\begin{enumerate}

\smallskip

\item $L_{\Delta_X}$ can be realized as an iterated cone of the branes $\fL_{\Delta_{X_\alpha},0}$.\

\smallskip

\item $\fL_{\Delta_{X_\alpha},1}$ is an external product of the form 
$L_{\alpha} \times\alpha_X(\mu_X(\CF_\alpha)).$\

\smallskip

\item 
For all $t\in [0,1]$, the functors
$$
\tilde\Psi_{\fL_{\Delta_{X_\alpha},t}*}|_{F(T^*X)_\Lambda}:{F(T^*X)_\Lambda} \to \rmod(F(T^*X))
$$
are quasi-isomorphic.

\end{enumerate}
The key point in explaining the third property is that the families $\fL_{\Delta_{X_\alpha},t}$
are {non-characteristic} with respect
to the conical Lagrangian $\Lambda\subset T^*X$ of the first factor.
Roughly speaking, for a given conical Lagrangian $\Lambda'\subset T^*X$,
we can arrange (after appropriate perturbations)
so that the boundaries at infinity of the families $f{\Delta_{X_\alpha},t}$ 
do not intersect the boundary of $\Lambda\times \Lambda'$.

With the preceding in hand, by applying the identity functor to any object of
$F(T^*X)_\Lambda$,
we immediately conclude that it is quasi-isomorphic to an object coming from $Sh_c(X)$.

\medskip

Independently of the above arguments, one can construct an explicit $A_\infty$-functor
$$
\pi_X: F(T^*X)\to Sh_c(X)
$$
which is  a quasi-inverse to $\mu_X$ (see Section~\ref{sect from branes to sheaves}
for details). 
To define $\pi_X$, for each open subset $U\hra X$, 
consider the corresponding costandard brane 
$
L_{U!} 
$ 
(see Section~\ref{sect costandard branes}).
The underlying Lagrangian of $L_{U!}$ can be taken to be
the graph $\Gamma_{-d\log m}$ for any non-negative function $m:X\to\R$
that vanishes precisely on the complement $X\setminus U$.

Given an object $L$ of $F(T^*X)$, the assignment
$$
U\mapsto \hom_{F(T^*X)}(L_{U!} \otimes \orient_X[-\dim X], L)
$$
defines 
a contravariant $A_\infty$-functor
from the category of open sets of $X$ to the dg category of chain complexes.
Without much difficulty, it is possible to reinterpret this as a constructible complex of sheaves 
on $X$ which we take to be $\pi_X(L)$. 
It follows quickly from the definitions that we have a canonical quasi-isomorphism of functors
$$
\pi_X\circ\mu_X\simeq \id_{\Sh_c(X)}
$$
confirming that $\pi_X$ and $\mu_X$ are quasi-inverses.
The construction of $\pi_X$ is very similar to some results of~\cite{KO1,KO2}
though presented in the language of constructible sheaves rather than F\'ary functors.

Finally, it is simple to understand basic properties of the functors  $\pi_X$ and $\mu_X$
such as their constructibility.
First, fix a stratification $\CS=\{S_\alpha\}$ of X, and consider the full subcategory
$$
Sh_\CS(X) \subset Sh_c(X)
$$
of complexes constructible with respect to $\CS$. Then by construction, we have
$$
\mu_X: 
Sh_\CS(X)\to F(T^*X)_{\Lambda_\CS}
$$ 
where $\Lambda_\CS= \cup_\alpha T^*_{S_\alpha} X$.
Conversely, given a conical Lagrangian $\Lambda\subset T^*X$, as part of the construction
of $\pi_X$, we verify that we have
$$
\pi_X:F(T^*X)_\Lambda\to Sh_\CS(X)
$$ 
for any stratification $\CS$ such that
$\Lambda\subset \Lambda_\CS$.
It is also simple to see that the functors  $\pi_X$ and $\mu_X$ 
interchange the brane duality $\alpha_X$ with Verdier duality $\CD_X$
and are also compatible with the basic operations on sheaves and branes.


\subsection{Applications}
We will postpone most applications to future papers and restrict ourselves here to one
immediate application to symplectic topology. 

The question of the possible structure of
compact Lagrangian submanifolds of $T^*X$ has seen some progress in recent years.
For a recent example of the subject, we refer the reader
to the paper of Seidel~\cite{Seid04}. It contains a brief summary of other relevant 
works of Lalonde-Sikarov~\cite{LS}, Viterbo~\cite{Vit}, and Buhovsky~\cite{Buh},
and is the paper in the subject which is closest to this one in methods. Namely, 
the main point is that we may reinterpret properties of objects of the Fukaya category of $T^*X$ 
in terms of the structure of their underlying Lagrangian submanifolds.

Consider a compact connected Lagrangian submanifold $L\subset T^*X$.
To ensure that we may
lift $L$ to an object
of the Fukaya category of $T^*X$, we assume first 
that $L$ is exact
and has trivial Maslov class. 
For simplicity in the following statement, we will also assume $\pi_1(X)$ is trivial. 
Thus in particular, the second Stiefel-Whitney class
$w_2(T^*X)$, which is the square of the pullback of 
$w_1(X)$, must vanish
so $T^*X$ is spin.
This also implies that $L$ is orientable
since the difference between $w_1(L)$ and the restriction of the pullback
of $w_1(X)$ is the 
$\Z/2\Z$-reduction of the Maslov class.
Further, we will assume that $w_2(L)$ is the restriction of the pullback of $w_2(X)$, 
so that $L$ is relatively spin with respect to the 
background class given by the pullback of $w_2(X)$. 

The following application of the equivalence of $F(T^*X)$ and $Sh_c(X)$ generalizes
part of the main statement of~\cite{Seid04} from the case when $X$ is a sphere.
During the preparation of this paper,
we learned that a similar characterization of compact exact branes has recently been obtained by 
Fukaya-Seidel-Smith~\cite{FSS, Smith} by a variety of different methods.

\begin{thm} Assume $X$ and $L$ are as above. Then $L$ is equivalent in the Fukaya
category of $T^*X$ to (a shift of) the zero section. From this, we conclude:
\begin{enumerate}
\item $[L] =  \pm [X] \in H_{\dim X}(T^*X,\C)$.

\item $H^*(L,\C) \simeq H^*(X,\C)$.

\item If $L'\subset T^*X$ is another Lagrangian submanifold satisfying the same conditions,
then we have a lower bound on the  (possibly infinite) number of intersection points
$$
\#(L\cap L' ) \geq \sum_k \dim H^k(L,\C).
$$

\end{enumerate}
\end{thm}

\begin{proof}
By assumption, we may equip $L$ with a brane structure so that it becomes
an object of $F(T^*X)$. Applying $\pi_X$ to this brane produces an object $\CF$ of $Sh_c(X)$.
Since $L$ is compact, $\CF$ is constructible with respect to the trivial stratification of $X$.
In other words, the cohomology sheaves of $\CF$ are local systems.
By assumption, $X$ is simply-connected, so these local systems are all trivial.

Applying $\pi_X$ to standard calculations in $F(T^*X)$, we see that 
$$
\Ext_{Sh_c(X)}^*(\CF,\CF)\simeq H^*(L).
$$ 
In particular, $H^m(L)=0$, for $m>\dim X$, implies $\Ext_{Sh_c(X)}^m(\CF,\CF)=0$, for $m>\dim X$.
By writing $\CF$ as a successive extension of (shifts of) the constant sheaf $\C_X$,
we see that this bound forces $\CF$ to reduce to (a shift of) $\C^{\oplus k}_X$ for some $k>0$.
But then since $L$ is connected, $H^0(L)\simeq \C$, and so $\CF$ is isomorphic to 
(a shift of) $\C_X$ itself.
Thus applying $\mu_X$, we conclude that $L$ is equivalent to (a shift of) the zero section 
in $F(T^*X)$. 
This implies assertions (2) and (3) immediately.

For assertion (1), by applying $\mu_X$ to the skyscraper sheaf $\C_{\{x\}}$,
we can consider the conormal to a point $T^*_{\{x\}} X$ as an object of $F(T^*X)$.
Since $\Ext_{Sh_c(X)}(\C_X,\C_{\{x\}})$ is isomorphic to $\C$ concentrated 
in degree zero, (after a possible shift)
the cohomology of $\hom_{F(T^*X)}(L,T^*_{\{x\}} X)$ is as well. In particular, the Euler characteristic
of the complex of intersection points is equal to $\pm 1$. This implies assertion (1).
\end{proof}


\subsection{Acknowledgements}
It is a pleasure to thank Eric Zaslow for many enlightening discussions,
and Paul Seidel for detailed remarks on the arguments.
I would also like to thank Sasha Beilinson and Yong-Geun Oh for comments on the exposition,  
and Kevin Costello, Ezra Getzler, and Katrin Wehrheim for very useful discussions.




\section{Constructible sheaves}\label{sec sheaves}

In this section, we first review background material on the constructible derived category,
then recall a differential graded model of it. Finally, we explain
how one can construct constructible sheaves
out of certain $A_\infty$-functors.


\subsection{Derived category}
In this section, we briefly recall the construction of the constructible derived category of
a real analytic manifold.
For a comprehensive treatment of this topic, the reader could consult
the book of Kashiwara-Schapira~\cite{KS}. 

\medskip
Let $X$ be a topological space. Let $\Top(X)$ be the category whose objects
are open sets $U\hra X$, and morphisms are 
inclusions $U_0 \hookrightarrow U_1$ of open sets:
$$
\hom_{\Top(X)} (U_0, U_1)
=
\left\{
\begin{array}{cl}
pt & \mbox{ when $U_0 \hra U_1$,} \\
\emptyset & \mbox{ when $U_0 \not\hra U_1$.}
\end{array}
\right.
$$

\medskip

Let $\Vect$ be the abelian category of complex vector spaces.

\medskip

The derived category of sheaves of complex vector spaces on $X$ is 
traditionally defined
via the following sequence of constructions:

\medskip
\noindent
1. {\em Presheaves}. Presheaves on $X$ are functors $\CF: \Top(X)^\circ\to \Vect$
where $\Top(X)^\circ$ denotes the opposite category. Given an open set $U\hra X$,
one writes $\CF(U)$ for the sections of $\CF$ over $U$, and given an inclusion $U_0 \hra U_1$
of open sets, one writes $\rho^{U_1}_{U_0}: \CF(U_1)\to \CF(U_0)$ for the
corresponding restriction map.

\medskip
\noindent
2. {\em Sheaves}. Sheaves on $X$ are presheaves $\CF: \Top(X)^\circ\to \Vect$
which are locally determined in the following sense.
For any open set $U\hra X$,
and covering $\fU=\{U_i\}$ of $U$ by open subsets $U_i\hra U$, there is a complex of vector spaces
$$
\xymatrix{
0 \ar[r] & \CF(U) \ar[r]^-{\delta} & \prod_i \CF(U_i) \ar[r]^-{\delta_0} & \prod_{i,j} \CF(U_i\cap U_j),
}
$$ 
where $\delta = \prod_i \rho^U_{U_i}$ and 
$\delta_0 = \prod_{i,j} \left (\rho^{U_i}_{U_i\cap U_j} - \rho^{U_j}_{U_i\cap U_j}\right )$.
A sheaf is a presheaf for which $\ker(\delta) = \ker(\delta_0)/\on{im}(\delta)=0$ for all
open sets 
and coverings of open sets.
%

\medskip

Sheaves on $X$ form an abelian category and thus one can continue with
the following sequence of general homological constructions:

\medskip
\noindent
3. {\em Complexes}. Let $C(X)$ be the abelian category of complexes of sheaves on $X$
with morphisms the degree zero chain maps.
Given a complex of sheaves $\CF$, one writes $H(\CF)$ for
the (graded) cohomology sheaf of $\CF$.

\medskip
\noindent
4. {\em Homotopy category}. Let $K(X)$ be the homotopy category of sheaves on $X$
with objects complexes of sheaves and morphisms homotopy classes of maps.
This is a triangulated category whose distinguished triangles are isomorphic to 
the standard mapping cones.

\medskip
\noindent
5. {\em Derived category}. The derived category $D(X)$ of sheaves on $X$
is defined to be the localization of $K(X)$ with respect to homotopy classes of quasi-isomorphisms (maps
inducing isomorphisms on cohomology).
Acyclic objects form a null system in $K(X)$, and thus $D(X)$ inherits the structure of 
triangulated category.

\medskip

With the derived category $D(X)$ in hand, one can define many variants by imposing
topological and homological conditions on objects.

\medskip
\noindent
6. {\em Bounded derived category}.
The bounded derived
category $D^\flat(X)$ is defined to be the full subcategory of $D(X)$ of bounded complexes.

Two standard equivalent descriptions are worth keeping in mind: first,
there is the more flexible description of $D^\flat(X)$ as
the full subcategory of $D(X)$ of complexes with bounded cohomology; second, 
there is the computationally useful description of
$D^\flat(X)$ as the homotopy category
of complexes
of injective sheaves with bounded cohomology.

\medskip
\noindent
7. {\em Constructibility}.  Assume $X$ is a real analytic manifold.
Fix an analytic-geometric category
$\CC$ in the sense of~\cite{vdDM}. For example, one could take $\CC(X)$ to be the
subanalytic subsets of $X$ as described in~\cite{BM}.

Let $\CS=\{S_\alpha\}$ be a Whitney stratification of $X$ by $\CC$-submanifolds
$i_\alpha: S_\alpha\hra X$.
An object $\CF$ of $D(X)$ is said to be $\CS$-constructible 
if the restrictions
$
i_\alpha^*H(\CF)
$
of its cohomology sheaf
to the strata of $\CS$ are finite-rank and locally constant.

The $\CS$-constructible derived category $D_\CS(X)$ is the full subcategory of $D(X)$ of $\CS$-constructible objects.
The constructible derived category $D_c(X)$ is the full subcategory of $D(X)$ of objects which are $\CS$-constructible for some Whitney stratification $\CS$.

Note that if the stratification $\CS$ is finite (for example, if $X$ is compact), then the finite-rank condition implies that all $\CS$-constructible objects have bounded cohomology. In other words,
within $D(X)$, 
every object of $D_\CS(X)$ is isomorphic to an object of $D^\flat(X)$.


\subsection{Differential graded category}
The derived category $D(X)$ is naturally the cohomology category of a differential
graded (dg) category $Sh(X)$. To define it, we will return to the sequence of 
homological constructions 
listed above and perform some modest changes. 
Two principles guide 
such definitions: (1) structures (such as morphisms and higher exts) should be defined at the level
of complexes not their cohomologies; 
and (2) properties (such as constructibility) should be imposed at the level of cohomologies rather than
complexes. The first principle ensures we will not lose important information,
while the second ensures we will have sufficient flexibility.
As an example of the latter, we prefer the realization of the bounded derived
category $D^\flat(X)$ as the full subcategory of $D(X)$ of complexes with bounded cohomologies
rather than of strictly bounded complexes.

The reader could consult~\cite{Keller, Drinfeld, KellerICM} 
for background on dg categories, in particular,
a discussion of the construction of dg quotients.

\medskip

Recall that sheaves on $X$ form an abelian category. The following sequence of
homological constructions can be performed on any abelian category:

\medskip
\noindent
1. {\em Dg category of complexes}. Let $C_{dg}(X)$ be the dg category
with objects complexes of sheaves and morphisms the usual complexes of maps between complexes.
In particular, the degree zero cycles in such a morphism complex are the usual degree zero
chain maps which are the morphisms of the ordinary category $C(X)$.

\medskip
\noindent
2. {\em Dg derived category}. The dg derived category $Sh(X)$
is defined to be the dg quotient of $C_{dg}(X)$ by the full subcategory of acyclic objects.
This is a triangulated dg category whose
cohomology category $H(Sh(X))$ is canonically equivalent (as a triangulated category) to the usual derived category $D(X)$.

\medskip

One can cut out full triangulated dg subcategories of $Sh(X)$ by specifying
full triangulated subcategories of its cohomology category $H(Sh(X))\simeq D(X)$. 

\medskip
\noindent
3. {\em Bounded dg derived category}.
The bounded dg
derived category $Sh^\flat(X)$ is defined to be the full dg subcategory of $Sh(X)$ 
of objects projecting to $D^\flat(X)$.

\medskip
\noindent
4.  {\em Constructibility}. 
 Assume $X$ is a real analytic manifold, and 
fix an analytic-geometric category $\CC$.
The constructible dg derived category $Sh_{c}(X)$ is the full dg subcategory of $Sh(X)$
of objects projecting
to $D_c(X)$. 
For a Whitney stratification $\CS$ of $X$, the $\CS$-constructible dg derived category 
$Sh_{\CS}(X)$ is the full dg subcategory of $Sh(X)$
of objects projecting
to $D_\CS(X)$.

\medskip

The formalism of Grothendieck's six (derived) operations $f^*, f_*, f_!, f^!, \intHom, \otimes$
can be lifted to the constructible dg derived category $Sh_c(X)$ 
(see for example~\cite{Drinfeld} for a general discussion of deriving
functors in the dg setting). In our case, one concrete approach
is to recognize that the natural map $C_{dg, c}(\mathfrak{Inj}(X)) \to Sh_c(X)$ 
from the dg category $C_{dg, c}(\mathfrak{Inj}(X))$ of complexes of injective sheaves
with constructible cohomology is a quasi-equivalence.
With this in hand, one can define derived functors by evaluating their
naive versions on $C_{dg, c}(\mathfrak{Inj}(X))$.
Since we will only consider derived functors, we will denote them by the above unadorned symbols.

\medskip

Throughout the remainder of this paper, we assume that $X$ is a real
analytic manifold,  and we
fix an analytic-geometric category $\CC$. All subsets will be $\CC$-subsets
unless otherwise stated.


\subsection{Standard bases}

We recall here several standard bases for the 
$\CT$-constructible dg derived category $Sh_\CT(X)$
for a triangulation $\CT=\{\tau_\fra\}$ of $X$ by simplices
$j_\fra:\tau_\fra\hra X$. (They are also well-known as basic examples 
in the theory of exceptional collections.)

\medskip

Define ${\CC_*}(\CT)$ to be the full dg category of $Sh_\CT(X)$ of
 standard objects
$j_{\fra *} \C_{\tau_\fra}$.
The morphisms between standard objects are quasi-isomorphic to
complexes concentrated in degree zero
$$
\hom_{Sh_\CT(X)} (j_{\frb *} \C_{\tau_\frb},j_{\fra *} \C_{\tau_\fra})
\simeq
\left\{
\begin{array}{cl}
\C & \mbox{ when $\tau_\fra \subset \ol\tau_\frb$,} \\
0 & \mbox{ when $\tau_\fra \not\subset \ol\tau_\frb$.}
\end{array}
\right.
$$
The composition maps are given by the linearization of the obvious poset relations.

\begin{lem}\label{standard envelope}
$Sh_\CT(X)$ is the triangulated envelope of ${\CC_*}(\CT)$.
\end{lem}

\begin{proof}

Let $i_{\geq k}:\CT_{\geq k}\hra X$ denote the union of the simplices of $\CT$
of dimension greater than or equal to $k$, and let $j_{< k}:\CT_{<k}\hra X$ denote the union of the simplices of $\CT$
of dimension less than $k$. 

 Let $Sh_{\CT_{\geq k} *}(X)$ denote the full dg
subcategory of $Sh_\CT(X)$ of objects of the form $\CF\simeq i_{\geq k*}\CF_{\geq k}$. 
By the standard triangle
$$
j_{< k!} j^!_{< k}\CF \to \CF\to i_{\geq k*} i_{\geq k}^*\CF\stackrel{[1]}{\to},
$$ 
this is 
equivalent to $j^!_{< k}\CF \simeq 0$. Let us show by induction that $Sh_{\CT_{\geq k} *}(X)$ is generated by the standard objects 
$j_{\fra *} \C_{\tau_\fra}$ associated to simplices $\tau_\fra$ with $\dim\tau_\fra \geq k$.
In particular, the assertion of the lemma is the case $k=0$.

For $k=\dim X$, $Sh_{\CT_{\geq k} *}(X)$ consists of complexes of standard
objects on the top-dimensional simplices and nothing else.

Suppose we know the assertion for all $\ell> k$.
For any object $\CF$ of $Sh_{\CT_{\geq k}*}(X)$,
we have a distinguished triangle
$$
 j_{<k+1*} j_{<k+1}^!\CF
\to \CF 
\to i_{\geq  {k+1}*} i_{\geq {k+1}}^*\CF
\stackrel{[1]}{\to}.
$$
By induction, 
we can express the object $i_{\geq  {k+1}*} i_{\geq {k+1}}^*\CF$
in terms of the standard objects associated to the simplices of $\CT_{\geq k+1}$.
Applying $j_{<k}^!$ to the above triangle, we obtain
$$j_{<k}^! j_{<k+1*} j_{<k+1}^!\CF 
\simeq 0,
$$
and so 
we can express the object $j_{<k+1*} j_{<k+1}^!\CF$
in terms of the standard objects associated to the simplices of dimension $k$.
\end{proof}


We can obtain a costandard basis by applying Verdier duality to the standard basis
as follows.

Define ${\CC_!}(\CT)$ to be the full dg category of $Sh_\CT(X)$ of
 costandard objects
$j_{\fra !} \omega_{\tau_\fra} \simeq \CD(j_{\fra *} \C_{\tau_\fra})$. 
Here $\omega_{\tau_\alpha}\simeq \CD(\C_{\tau_\alpha})$ denotes the Verdier dualizing
complex of $\tau_\alpha$; it is canonically isomorphic to the shifted orientation sheaf
$\orient_{\tau_\alpha}[\dim\tau_\alpha]$, and so, given a choice of orientation of $\tau_\alpha$,
isomorphic to the 
shifted constant sheaf $\C_{\tau_\alpha}[\dim\tau_\alpha]$.
The morphisms between costandard objects are quasi-isomorphic to
complexes concentrated in degree zero
$$
\hom_{Sh_\CT(X)} (j_{\fra !} \omega_{\tau_\fra},j_{\frb !} \omega_{\tau_\frb})
\simeq
\left\{
\begin{array}{cl}
\C & \mbox{ when $\tau_\fra \subset \ol\tau_\frb$,} \\
0 & \mbox{ when $\tau_\fra \not\subset \ol\tau_\frb$.}
\end{array}
\right.
$$
The composition maps are given by the linearization of the obvious poset relations.

Since Verdier duality is an anti-equivalence, Lemma~\ref{standard envelope} implies
the following.

\begin{lem}\label{costandard envelope}
$Sh_\CT(X)$ is the triangulated envelope of ${\CC_!}(\CT)$.
\end{lem}


Here is an alternative basis of standard objects associated to open sets.
For each simplex $\tau_\fra$ of $\CT$, let $i_\fra:st_\fra\hra X$ be its star 
$$
st_\fra = \bigsqcup_{ \tau_\fra \subset\ol\tau_\frb} \tau_\frb.
$$
Note that $st_\fra$ is an open contractible submanifold of $X$,
 and we have $\tau_\fra\subset \ol\tau_\frb$
if and only if
$st_\frb\subset st_\fra$.

Define $\CC_{st*}(\CT)$ to be the full dg category of $Sh_\CT(X)$ of
 standard objects
$i_{\fra *} \C_{st_\fra}$.
The morphisms between standard objects are quasi-isomorphic to
complexes concentrated in degree zero
$$
\hom_{Sh_\CT(X)} (i_{\fra *} \C_{st_\fra},i_{\frb *} \C_{st_\frb})
\simeq
\left\{
\begin{array}{cl}
\C & \mbox{ when $\tau_\fra \subset \ol\tau_\frb$,} \\
0 & \mbox{ when $\tau_\fra \not\subset \ol\tau_\frb$.}
\end{array}
\right.
$$
The composition maps are given by the linearization of the obvious poset relations.

\begin{lem}\label{standard star envelope}
$Sh_\CT(X)$ is the triangulated envelope of ${\CC}_{st*}(\CT)$.
\end{lem}

\begin{proof}
We continue with the notation of the proof of Lemma~\ref{standard envelope}.

Let us show by induction that $Sh_{\CT_{\geq k} *}(X)$ is generated by the standard objects 
$i_{\fra *} \C_{st_\fra}$ associated to the stars of simplices with $\dim\tau_\fra \geq k$.
Recall that the proof of Lemma~\ref{standard envelope} shows that
$Sh_{\CT_{\geq k} *}(X)$ is generated by the standard objects $j_{\fra *} \C_{\tau_\fra}$
associated to simplices with $\dim\tau_\fra \geq k$.

For $k=\dim X$,
if $\dim\tau_\fra =\dim X$, then $st_\fra = \tau_\fra$ and so 
$j_{\fra *} \C_{\tau_\fra} = i_{\fra *} \C_{st_\fra}$.

Suppose we know the assertion for all $k>\dim\tau_\fra$. Let $i'_\fra:st'_\fra\hra X$
be the punctured star $st'_\fra = st_\fra \setminus \tau_\fra$,
and consider the distinguished triangle
$$
j_{\fra *}\nu_{\tau_\fra}
\to i_{\fra*} \C_{st_\fra}\to 
i'_{\fra*}\C_{st'_\fra}
\stackrel{[1]}{\to}.
$$
Here $\nu_{\tau_\alpha}\simeq j^!_\fra\C_X$ 
is canonically isomorphic to the shifted normal orientation sheaf
$\orient_{X/\tau_\alpha}[-\dim\tau_\alpha]$, and so, given a choice of 
normal orientation,
isomorphic to the 
shifted constant sheaf $\C_{\tau_\alpha}[-\dim\tau_\alpha]$.

By induction, we can express $i'_{\fra*}\C_{st'_\fra}$ in terms of
 the standard objects 
$i_{\frb *} \C_{st_\frb}$ associated to the stars of simplices with $\dim\tau_\frb \geq k$.
Therefore we can express the standard object $j_{\fra *}\C_{\tau_\fra}$ as well.
\end{proof}


We can obtain a costandard basis by applying Verdier duality to the above standard basis
as follows.

Define $\CC_{st!}(\CT)$ to be the full dg category of $Sh_\CT(X)$ of
 costandard objects
$i_{\fra !} \omega_{st_\fra} \simeq \CD(i_{\fra *} \C_{st_\fra})$.
The morphisms between costandard objects are quasi-isomorphic to
complexes concentrated in degree zero
$$
\hom_{Sh_\CT(X)} (i_{\frb !} \omega_{st_\frb},i_{\fra !} \omega_{st_\fra})
\simeq
\left\{
\begin{array}{cl}
\C & \mbox{ when $\tau_\fra \subset \ol\tau_\frb$,} \\
0 & \mbox{ when $\tau_\fra \not\subset \ol\tau_\frb$.}
\end{array}
\right.
$$
The composition maps are given by the linearization of the obvious poset relations.

Since Verdier duality is an anti-equivalence, Lemma~\ref{standard star envelope} implies
the following.

\begin{lem}\label{costandard star envelope}
$Sh_\CT(X)$ is the triangulated envelope of ${\CC}_{st!}(\CT)$.
\end{lem}


\subsection{Quasi-representable modules}\label{representable section}
This section is not logically needed in what follows.
We include it to give a feeling for 
the information contained in an $A_\infty$-module over 
the constructible dg derived category $Sh_c(X)$.
By restricting a quasi-representable $A_\infty$-module over $Sh_c(X)$
to the costandard sheaves of open subsets, one obtains what could be called
an $A_\infty$-sheaf.
For background on $A_\infty$-categories, the reader could consult~\cite{Seidel}
and the references therein.
Note that any dg category can be viewed as an an $A_\infty$-category
with vanishing higher compositions.

\medskip

Let $\Ch=C_{dg}(pt)$ denote the dg category of chain complexes of complex vector spaces,
and let
$ \rmod(Sh_{c}(X))$ denote the $A_\infty$-category of $A_\infty$-functors
$$
\CM:Sh_c(X)^\circ\to \Ch.
$$
In keeping with usual nomenclature (inspired by considering categories with a single object), 
we will also refer to such functors as right $Sh_c(X)$-modules. 
The Yoneda functor 
$$
\CY_{r}:Sh_{c}(X) \to \rmod(Sh_{c}(X))
$$
$$
\CY_{r}(\CF) : \CP \mapsto \hom_{Sh_c(X)}(\CP,\CF)
$$
is a quasi-embedding of $A_\infty$-categories in the sense that
the induced cohomology functor
is fully faithful.
We will say that a module $\CM$ is quasi-represented by
an object $\CF$ if there is a quasi-isomorphism of modules
$\CY_r(\CF)\simeq \CM$.

\medskip

Fix a Whitney stratification $\CS=\{S_\alpha\}$ of $X$ by submanifolds
$i_\alpha:S_\alpha\hra X$.
We would like to characterize when a right $Sh_c(X)$-module
is quasi-represented by an object of $Sh_\CS(X)$. 
The following two properties clearly hold for any right $Sh_c(X)$-module $\CM$
quasi-represented by an object of $Sh_\CS(X)$.

\medskip
\noindent
({\em f-r}) 
For any submanifold $i:Y\hra X$ with costandard object $i_{!} \C_Y$,
the cohomology of the evaluation 
$\CM(i_{!} \C_Y) 
$ 
is finite-rank.

\medskip

\medskip
\noindent
($\CS${\em-lc}) Let $\CT=\{\tau_\fra\}$ be any triangulation of $X$ by simplices 
$j_\fra: \tau_\fra \hra X$ refining the stratification $\CS$
in the sense that each stratum of $\CS$ is a union of simplices of $\CT$.

For each pair $\tau_\fra \subset \ol \tau_\frb$ of simplices of $\CT$ 
such that $\tau_\fra,\tau_\frb$ both lie in a single stratum
of $\CS$, 
the natural evaluation map
$$
\xymatrix{
\CM( i_{\frb!} \omega_{\tau_\frb}) 
\ar[r] 
&
\CM(i_{\fra !} \omega_{\tau_\fra})
}
$$
is a quasi-isomorphism.

\medskip

Clearly any object $\CF$ of $Sh_c(X)$ whose Yoneda module $\CY_r(\CF)$
satisfies ($\CS${\em -lc}) belongs to $Sh_\CS(X)$.
The remainder of this section is devoted to showing that any right $Sh_c(X)$-module
satisfying ({\em f-r}) 
and ($\CS${\em -lc}) is quasi-represented by an object of $Sh_\CS(X)$.
%


\begin{lem}\label{representable for triangulation}
For any right $Sh_c(X)$-module $\CM$ satisfying the finite-rank condition (f-r),
and any triangulation $\CT$,
there is an object $\CF_\CT$ of $Sh_\CT(X)$
that quasi-represents the restriction
$\CM|_{Sh_\CT(X)}$.
\end{lem}

\begin{proof}
Recall the full dg category $\CC_!(\CT) \subset Sh_c(X)$ 
of costandard objects introduced in the preceding section.
By Lemma~\ref{costandard envelope}, the embedding realizes $Sh_\CT(X)$ as the triangulated envelope
of ${\CC_!}(\CT)$, and hence
the corresponding restriction of modules
is a quasi-embedding as well
$$
\rmod(Sh_\CT(X)) \hra \rmod({\CC_!}(\CT)).
$$
Thus to prove the lemma, it suffices to define an object $\CF_\CT$ of $\Sh_\CT(X)$
along with a quasi-isomorphism of modules
$$
\CY_r(\CF_\CT)|_{{\CC_!}(\CT)}\simeq \CM|_{{\CC_!}(\CT)}
$$

Consider the (oriented) standard objects
$j_{\fra *} \omega_{\tau_\fra}$ associated to the simplices $j_\fra:\tau_\fra\hra X$.
They provide the simplest right ${\CC_!}(\CT)$-modules in the sense that
$$
\hom_{\Sh_\CT(X)}( j_{\frb !} \omega_{\tau_\frb},  j_{\fra *} \omega_{\tau_\fra})\simeq
\hom_{\Sh_\CT(X)}( \omega_{\tau_\frb},  j^!_\frb j_{\fra *} \omega_{\tau_\fra})\simeq
\left\{
\begin{array}{ll}
\C & \mbox{ when $\fra=\frb$,} \\
0 & \mbox{ when $\fra \not = \frb$.}
\end{array}
\right.
$$
We will use this to show that $\CM|_{{\CC_!}(\CT)}$ can be expressed as 
an iterated cone of shifts of the standard modules 
$\CY_r(j_{\fra *} \omega_{\tau_\fra})|_{{\CC_!}(\CT)}$.

The argument, similar to that of Lemma~\ref{standard envelope},
 is an induction on the dimension of the simplices of $\CT$, beginning with the open simplices.
Let $n=\dim X$, and for $0\leq k\leq n$, let $\frI_k$ be those indices labelling 
simplices of $\CT$ of dimension $k$.

\medskip
\noindent
({\em Step} $n$) For $\fra_n\in\frI_n$, consider the evaluation 
$\CM({\fra_n}) = \CM(j_{\tau_{\fra_n }!}\omega_{\tau_{\fra_n}})$.

We claim that there is a canonical map of right ${\CC_!}(\CT)$-modules
$$
q_n:\CM\to \sum_{\fra_n\in\frI_n} \CY_r(j_{{\fra_n}*} (\CM({\fra_n})\otimes\omega_{\tau_{\fra_n}})).
$$
To see this, note that for any $\tau_\frb$, the right hand side evaluates to be
$$
\left(
\sum_{\fra_n\in\frI_n} \CY_r(j_{{\fra_n}*} (\CM({\fra_n})\otimes\omega_{\tau_{\fra_n}}))
\right)
(j_{{\frb} !}\omega_{\tau_{\frb}})\simeq
\left\{
\begin{array}{ll}
\CM({\fra_n}) & \mbox{ when $\frb = \fra_n$,} \\
0 & \mbox{ when $\frb \not = \fra_n$.}
\end{array}
\right.
$$
Thus we can take the first-order part of $q_n$ to be the identity when $\frb=\fra_n$ and zero otherwise.
Furthermore, there are no higher-order terms to define.

Let $\CM_{<n}$ denote the cone of $q_n$. By construction, since $q_n$ evaluated at any 
$j_{\fra !}\omega_{\tau_\fra}$, for any
$\fra\in\frI_n$, is
an isomorphism at the chain level, we can arrange so that $\CM_{<n}(j_{\fra !}\omega_{\tau_\fra}) = 0$, 
for all $\fra\in\frI_n$.

\medskip
\noindent
({\em Step} $n-1$) For $\fra_{n-1}\in\frI_{n-1}$, consider the evaluation $\CM_{<n}(\fra_{n-1}) = 
\CM_{<n}(j_{\fra_{n-1} !}\omega_{\tau_{\fra_{n-1}}})$.
We claim that there is a canonical map of right ${\CC_!}(\CT)$-modules
$$
q_{n-1}:\CM_{<n}\to \sum_{\fra_{n-1}\in\frI_{n-1}} \CY_r(j_{\fra_{n-1}*} (\CM_{<n}(\fra_{n-1})
\otimes\omega_{\tau_{\fra_{n-1}}})).
$$
To see this, note that for any $\tau_\frb$, the left hand side evaluates to be
$$
\left(
\sum_{\fra_{n-1}\in\frI_{n-1}} \CY_r(j_{\fra_{n-1}*} (\CM({\fra_{n-1}})\otimes\omega_{\tau_{\fra_{n-1}}}))
\right)
(j_{\frb !}\omega_{\tau_\frb})\simeq
\left\{
\begin{array}{ll}
\CM_{<n}({\fra_{n-1}}) & \mbox{ when $\frb = \fra_{n-1}$,} \\
0 & \mbox{ when $\frb \not = \fra_{n-1}$.}
\end{array}
\right.
$$
Thus we can take the first-order part of $q_{n-1}$ to be the identity when $\frb=\fra_{n-1}$ 
and zero otherwise.
Furthermore, there are no higher-order terms to define.

Let $\CM_{<n-1}$ denote the cone of $q_{n-1}$. By construction, since $q_{n-1}$ 
evaluated at 
$j_{\fra !}\omega_{\tau_\fra}$, for any $\fra\in\frI_{n-1}\cup \frI_n$, is
an isomorphism at the chain level, we can arrange so that 
$\CM_{<n-1}(j_{\fra !}\omega_{\tau_\fra}) = 0$, for all $\fra\in\frI_{n-1}\cup \frI_n$.

\medskip

And so on. In the end, we see that $\CM$ can be expressed by a finite sequence of 
cones of Yoneda modules of shifted standard objects.
\end{proof}

Consider a second  triangulation $\CT'$ of $X$ that refines $\CT$
in the sense that each simplex of $\CT$ is a union of simplices of $\CT'$.
So we have fully faithful dg embeddings
$$
Sh_\CT(X)\hra Sh_{\CT'}(X)\hra Sh_c(X),
$$
and the corresponding restriction of modules
$$
\rmod(Sh_c(X)) \to \rmod(Sh_{\CT'}(X))\to \rmod(Sh_{\CT}(X)).
$$

\begin{lem}\label{compatibility of representing objects}
Suppose the right $Sh_c(X)$-module $\CM$ also satisfies the $\CS$-locally constant
property ($\CS${-lc}) (in addition to (f-r)), and that the triangulation $\CT$ refines the stratification $\CS$.
Then for any triangulation $\CT'$ that refines $\CT$, there is a canonical isomorphism 
$\CF_\CT \simeq \CF_{\CT'}$ of objects of $Sh_c(X)$.
\end{lem}

\begin{proof}
Recall that for any right $Sh_c(X)$-module $\CM$ satisfying the finite-rank condition ({\em f-r}),
and any triangulation $\CT$,
Lemma~\ref{representable for triangulation} provides
an object $\CF_\CT$
along with a quasi-isomorphism of right ${Sh_\CT(X)}$-modules
$$
\CY_r(\CF_\CT)
\simeq 
\CM|_{Sh_\CT(X)}.
$$

Applying Lemma~\ref{representable for triangulation} to $\CT'$ 
and restricting to $Sh_\CT(X)$ provides
a quasi-isomorphism of right ${Sh_\CT(X)}$-modules
$$
\CY_r(\CF_{\CT'})|_{Sh_\CT(X)}
\simeq 
\CM|_{Sh_\CT(X)}.
$$

Finally, composing the above quasi-isomorphisms gives a
quasi-isomorphism of right ${Sh_\CT(X)}$-modules
$$
\CY_r(\CF_{\CT'})|_{Sh_\CT(X)}
\simeq 
\CY_r(\CF_{\CT}).
$$

With the above quasi-isomorphism in hand, to prove the lemma,
 it suffices to show that $\CF_{\CT'}$ is in fact an object of $Sh_\CT(X)$.
For each pair $\tau'_\fra \subset \ol \tau'_\frb$ of simplices of $\CT'$ 
such that $\tau'_\fra,\tau'_\frb$ both lie in a single stratum
of $\CS$, 
%
there is a diagram which commutes at the level of cohomology
$$
\xymatrix{
\CM( i_{\frb!} \omega_{\tau_\frb}) 
\ar[r] \ar[d]_-{\wr}
& \ar[d]^-{\wr}
\CM(i_{\fra !} \omega_{\tau_\fra})\\
\hom_{Sh_c(X)}( i_{\frb!} \omega_{\tau_\frb}, \CF_{\CT'}) \ar[r] & 
\hom_{Sh_c(X)}( i_{\fra!} \omega_{\tau_\fra}, \CF_{\CT'}) }
$$
The assumption ($\CS${\em -lc}) on $\CM$ implies that the upper horizontal arrow
is a quasi-isomorphism, and so the lower horizontal arrow is as well.
Thus the object $\CF_{\CT'}$
is in fact $\CT$-constructible.
\end{proof}

We summarize the preceding development in the following statement.

\begin{prop}\label{representability of module}
Let $\CM$ be an object of $\rmod(Sh_c(X))$ satisfying the properties ({f-r}) 
and ($\CS${-lc}). Then $\CM$ is quasi-represented by an object of $Sh_\CS(X)$.
\end{prop}

\begin{proof}
Fix a triangulation $\CT$ refining the stratification $\CS$. 
By Lemma~\ref{representable for triangulation}, there is an object
$\CF_\CT$ quasi-representing the restriction $\CM|_{Sh_\CT(X)}$.
Given any stratification $\CS'$, we can find a triangulation $\CT'$
that simultaneously refines $\CT$ and $\CS'$.
Thus by Lemma~\ref{compatibility of representing objects}, the object $\CF_\CT$
quasi-represents $\CM$ on all of $Sh_c(X)$. Finally, $\CM$ satisfies
property ($\CS${\em-lc}) and hence the Yoneda $Sh_c(X)$-module $\CY_r(\CF_\CT)$ does as well. This clearly
implies
that $\CF_\CT$ is $\CS$-constructible.
\end{proof}



\section{Microlocal branes}\label{sec fukaya}
In the sections below, we review some basic aspects of the Fukaya category of the cotangent bundle $T^*X$ from \cite{NZ}.
In particular, we discuss how constructible sheaves on $X$
embed into its triangulated envelope. 
We also collect some technical results on isotopies and standard branes
needed in what follows.


\subsection{Preliminaries}\label{sect prels}

In what follows, we work with a fixed compact real analytic manifold $X$
with cotangent bundle $\pi:T^*X\to X$. We often denote points of $T^*X$ by pairs
$(x,\xi)$ where $x\in X$ and $\xi\in T^*_x X$.
The material of this section is a condensed version
of the discussion of~\cite{NZ}. 

\medskip

Let $\theta\in\Omega^1(T^*X)$ denote the canonical one-form $\theta(v) =\xi(\pi_*v)$, for
$v\in T_{(x,\xi)}(T^*X)$,  and let 
$\omega =d\theta\in\Omega^2(T^*X)$ denote the 
canonical symplectic structure. 
For a fixed Riemannian metric on $X$, let $|\xi|:T^*X\to \R$ denote
the corresponding fiberwise linear length function.


\subsubsection{Compactification} To better control noncompact Lagrangians
in $T^*X$, it is useful to work with 
the cospherical compactification $\ol\pi:\ol T^*X\to X$ of the projection
$\pi:T^*X\to X$ obtained by attaching the cosphere bundle at infinity $\pi^\infty:T^\infty X\to X$.

Concretely, we can realize the compactification $\ol T^*X$ as the quotient
$$
\ol T^*X = \left((T^* X \times \R_{\geq 0})\setminus (X \times \{0\})\right)/ \R_+$$
where
$\R_{+}$ acts 
by dilations on both factors.
The canonical inclusion $T^*X\hookrightarrow \ol T^* X$ sends
a covector $\xi$ to the class of $[\xi,1]$.
The boundary at infinity 
$T^\infty X=\ol T^*X\setminus T^*X
$ 
consists of classes of the form $[\xi, 0]$ with $\xi$ a non-zero covector.
Given a Riemannian metric on $X$, one can identify
$\ol T^*X$ with the closed unit disk bundle $D^*X$,
and $T^\infty X$ with the unit cosphere bundle  $S^* X$,
via the map
$$[\xi,r]\mapsto (\hat \xi, \hat r),
\mbox{ where }|\hat\xi|^2 + \hat r^2=1.
$$

The boundary at infinity $T^\infty X$ carries a canonical contact 
distribution $\kappa \subset T (T^\infty X)$ with a well-defined notion of positive normal direction.
Given a Riemannian metric on $X$, under the induced identification of $T^\infty X$
with the unit cosphere bundle $S^* X$, the distribution $\kappa$ is
the kernel of the restriction of $\theta$.


\subsubsection{Conical almost complex structure}
\label{warp}
To better control holomorphic disks in $T^*X$, it is useful to work with 
an almost complex structure $J_{con}\in\on{End}(T(T^*X))$ 
which near infinity is invariant under dilations.

A fixed Riemannian metric on $X$ provides a canonical splitting
$
T(T^*X) \simeq T_{b} \oplus T_{f},
$
where $T_b$ denotes the horizontal base directions and $T_f$ the 
vertical fiber directions,
along with a canonical isomorphism 
$
j_0:T_b {\to} T_f
$
of vector bundles over $T^*X$.
We refer to the resulting almost complex structure
$$
J_{Sas}=
\left(
\begin{matrix}
0 & j_0^{-1} \\
-j_0 & 0
\end{matrix}
\right)
\in\on{End}(T_b \oplus T_f).
$$
as the Sasaki almost complex structure,
since by construction, the Sasaki metric on $T^*X$
is given by $g_{Sas}(v,v) = \omega(u, J_{Sas} v)$.

Fix positive constants $r_0, r_1>0$,
a bump function $b: \R\to \R$ such that $b(r) = 0$ for $r <r_0$, and $b(r) = 1$, for $r>r_1$,
and
set 
$w(x,\xi) = |\xi|^{b(|\xi|)},
$ where as usual
$|\xi|$ denotes the length of a covector with respect to the original
metric on $X$. 
We refer to the compatible almost complex structure 
$$
J_w=
\left(
\begin{matrix}
0 & w^{-1} j_0^{-1} \\
-w j_0 & 0
\end{matrix}
\right)
\in\on{End}(T_b \oplus T_f)
$$
as a(n asymptotically) conical almost complex structure
since near infinity $J_{con}$ is invariant under dilations.
The corresponding metric $g_{con}(u,v) = \omega(v, J_{con} v)$ presents
$T^*X$ near infinity as a metric cone over the unit {cosphere} bundle $S^*X$
equipped with the Sasaki metric.

One can view the conical metric $g_{con}$ as being compatible with the compactification $\ol T^*X$
in the sense that near infinity
it treats base and angular fiber directions on equal footing. Near infinity the metrics
on the level sets of $|\xi|$ are given by scaling the Sasaki metric
on the unit cophere bundle by the factor $|\xi|^{1/2}$.


\subsection{Brane structures}\label{sec branes}
By a Lagrangian $j:L\hookrightarrow T^*X$, we mean a closed (but not necessarily
compact) half-dimensional submanifold such that $TL$ is isotropic for the symplectic form $\omega$.
One says that $L$ is exact if the pullback of the one-form
$j^*\theta$ is cohomologous to zero. 

\medskip

By a brane structure on a Lagrangian ${L} \hra{T}^*X$, we mean a three-tuple
$
(\CE, \tilde\alpha,\flat)
$
consisting of a flat (finite-dimensional) vector bundle $\CE\to L$,
along with a grading $\tilde\alpha:L\to \R$ (with respect to the canonical
bicanonical trivialization) and a relative pin structure $\flat$ (with respect to the background class
$\pi^*(w_2(X))$. To remind the interested reader,
we include below a short summary of what the latter two structures entail.


\subsubsection{Gradings}
The almost complex structure $J_{con}\in\on{End}(T(T^*X))$ provides a holomorphic canonical bundle
$
\kappa= (\wedge^{\dim X}T^{hol}(T^*X))^{-1}.
$
According to~\cite{NZ}, there is a canonical trivialization $\eta^2$ of
the bicanonical bundle $\kappa^{\otimes 2}$ (and a canonical trivialization of $\kappa$
itself
if $X$ is assumed oriented).
Consider the bundle of Lagrangian planes ${\mathcal Lag}_{T^*X}\to T^*X$,
and the squared phase map
$$
\alpha : {\mathcal Lag}_{T^*X} \rightarrow U(1)
$$
$$
\alpha(\mathcal L)= \eta(\wedge^{\dim X} \mathcal L)^2/ |\eta(\wedge^{\dim X} \mathcal L)|^2.
$$

For a Lagrangian $L\hra  T^*X$ and a point $x\in L$, we obtain a map $\alpha:L\to U(1)$
by setting $\alpha(x)=\alpha(T_xL).$
The Maslov class $\mu(L)\in H^1(L,\Z)$
is the obstruction class
$
\mu = \alpha^*(dt),
$ 
where $dt$ denotes the standard one-form on $U(1)$.
Thus $\alpha$ has a lift to a map $\widetilde\alpha:L\to \R$
if and only if $\mu = 0$, and choices of a lift form a torsor over the group $H^0(L,\Z)$.
Such a lift $\widetilde \alpha:L\to \R$ is called a {grading} of the Lagrangian $L\hra T^*X$.


\subsubsection{Relative pin structures}
Recall that the group $Pin^+(n)$ is the double cover of
$O(n)$ with center $\Z/2\Z \times \Z/2\Z.$
A pin structure on a Riemannian manifold $L$ is a lift
of the structure group of $TL$ to $Pin^+(n).$
The obstruction to a pin structure is the second
Stiefel-Whitney class $w_2(L)\in H^2(L,\Z/2\Z)$,
and choices of pin structures form a torsor over the group
$H^1(L,\Z/2\Z)$. 

A relative pin structure on a
submanifold $L \hookrightarrow M$ with background class $[w]\in H^2(M,\Z/2\Z)$
can be defined as follows.  Fix a {\v C}ech cocycle $w$ representing $[w]$,
and let $w|_L$ be its restriction to $L$. Then a pin structure on $L$ relative to $[w]$ can be
defined to be  
an $w|_L$-twisted pin structure on $TL$. Concretely, this can be represented by a $Pin^+(n)$-valued
\v Cech $1$-cochain on $L$ whose coboundary is $w|_L$. 
Such structures are canonically independent of the choice of \v Cech representatives.

For Lagrangians $L\hra T^*X$, we will always consider relative pin structures $\flat$ on $L$
with respect to the fixed background class $\pi^*(w_2(X))\in H^2(T^*X,\Z/2\Z)$.


\subsection{Fukaya category}\label{sect Fukaya category}
We recall here the construction of the Fukaya $A_\infty$-category 
of the cotangent bundle $T^*X$. Our aim is not to review all of the details, but only those relevant to our later proofs.
For more details, the reader could consult~\cite{NZ} and the references therein.
In technical terms, the construction is a close relative of the category of vanishing cycles 
proposed by Kontsevich~\cite{Kontlecture} and Hori-Iqbal-Vafa~\cite{HIV},
and developed by Seidel~\cite{SeidmutationI}, \cite{SeidmutationII}, \cite{Seidel}.


\subsubsection{Objects}

An object of the Fukaya category of $T^*X$ is a four-tuple
$
(L, \CE, \tilde\alpha,\flat)
$
consisting of an exact (not necessarily compact) closed Lagrangian submanifold ${L} \hra{T}^*X$
equipped with a brane structure:
this includes a flat vector bundle $\CE\to L$,
along with a grading $\tilde\alpha:L\to \R$ 
(with respect to the canonical
bicanonical trivialization) and a relative pin structure $\flat$ (with respect to the background class
$\pi^*(w_2(X))$. 

To ensure reasonable behavior near infinity, we place two assumptions on the
Lagrangian $L$.
First, 
consider the compactification $\ol T^*X$ obtained by adding to $T^*X$ the cosphere bundle
at infinity $T^\infty X$. 
Then we fix an
analytic-geometric category $\CC$ once and for all, and
assume that
the closure $\ol L \hra  \ol T^*X$ is a $\CC$-subset. Along with other nice properties, 
this implies the following two key facts: 
\begin{enumerate}
\item The boundary at infinity 
$$
L^\infty =\ol L\cap T^\infty X
$$
is an isotropic subset of $T^\infty X$ with respect to the induced contact structure.
\item There is a real number $r>0$
such that the restriction of the length function 
$$|\xi| :L \cap \{|\xi| > r\}\to \R$$
has no critical points.
\end{enumerate}
As discussed below, the above properties guarantee we can make sense of ``intersections at infinity".

Second, to have a manageable theory of pseudoholomorphic maps with boundary on
such Lagrangians,
we also assume the existence of a perturbation $\psi$ that moves the initial
Lagrangian $L$ to a nearby Lagrangian tame (in the sense of~\cite{Sikorav})
with respect to the conical metric $g_{con}$. As confirmed in the Appendix,
all such perturbations lead to equivalent calculations.

We use the term Lagrangian brane to refer to objects of the Fukaya category.
When there is no chance for confusion,
we often write $L$ alone to signify the Lagrangian brane.


\subsubsection{Morphisms}

To define the morphisms between two branes, we must perturb Lagrangians
so that their intersections occur in some bounded domain.
To organize the perturbations, we recall the inductive
notion of a fringed set $R_{d+1}\subset \R^{d+1}_+$.
A fringed set $R_1\subset\R_+$ is any interval of the form $(0,r)$ for some $r>0$. 
A fringed set $R_{d+1}\subset\R_+^{d+1}$ is a subset satisfying the following:
\begin{enumerate}
\item $R_{d+1}$ is open in $\R^{d+1}_+$.
\item Under the projection $\pi:\R^{d+1}\to \R^d$ forgetting the last coordinate, 
the image $\pi(R_{d+1})$ is a fringed set.
\item If $(r_1,\ldots, r_d, r_{d+1})\in R_{d+1}$, then $(r_1,\ldots, r_d, r'_{d+1})\in R_{d+1}$
for $0<r'_{d+1}< r_{d+1}$.
\end{enumerate}

A Hamiltonian function $H:T^*X\to \R$ is said to be {controlled} if there is a real number
$r>0$ such that in the region $|\xi|>r$ we have $H(x,\xi)=|\xi|$. 
The corresponding Hamiltonian isotopy 
$\varphi_{H,t}:T^*X\to T^*X$ equals the normalized geodesic flow $\gamma_t$
in the region $|\xi|>r$.

As explained in~\cite{NZ},
given Lagrangians branes $L_0,\ldots, L_d\subset T^*X$, 
and controlled Hamiltonian functions $H_0,\ldots, H_d$,
we may choose a fringed set $R\subset \R^{d+1}$
such that for $(\delta_d,\ldots, \delta_0)\in R$, there is a real number $r>0$
such that for any $i\not = j$, we have
$$
\varphi_{H_i,\delta_i}(\ol L_i)\cap 
\varphi_{H_j,\delta_j}(\ol L_j) 
\quad
\mbox{lies in the region $|\xi|<r$.}
$$
By a further compactly supported Hamiltonian perturbation,
we may also arrange so that the intersections are transverse.

We consider finite collections of Lagrangian branes $L_0,\ldots, L_d\subset T^*X$ 
to come equipped with such perturbation data,
with the brane structures $(\CE_i, \tilde\alpha_i, \flat_i)$ and taming perturbations $\psi_i$ 
transported via the perturbations. Note that the latter makes sense since the normalized
geodesic flow $\gamma_t$ is an isometry of the metric $g_{con}$.
Then for branes $L_i,L_j$ with $i<j$,
the graded vector space of morphisms between them is defined to be
$$
{hom}_{F(T^*X)}(L_i,L_j) = \bigoplus_{p\in \psi_i(\varphi_{H_i,\delta_i}(L_i))\cap 
\psi_j(\varphi_{H_j,\delta_j}(L_j)) }
{\mathcal Hom}(\mathcal E_i \vert_{p},\mathcal E_j\vert_{p})[-\deg(p)].
$$
where the integer
$\deg(p)$ denotes the Maslov grading of the linear
Lagrangian subspaces at the
intersection.

It is worth emphasizing that near infinity the salient aspect of the above perturbation procedure 
is the relative position of the perturbed branes rather than their absolute position.
The following informal viewpoint
can be a useful mnemonic to keep the conventions straight.
In general, we always think of morphisms as ``propagating forward in time".
Thus 
to calculate the morphisms $\hom_{F(T^*X)}(L_0, L_1)$,
we have required that $L_0,L_1$ are perturbed near infinity by normalized geodesic flow
so that $L_1$ is further in the future than $L_0$. 
But what is important is not that they are both perturbed forward in time,
only that $L_1$ is further along the timeline than $L_0$. So for example,
we could perturb $L_0,L_1$ near infinity by normalized anti-geodesic flow
as long as $L_0$ is further in the past than $L_1$.


\subsubsection{Compositions}

Signed counts of
pseudoholomorphic polygons provide the differential and higher composition maps
of the $A_\infty$-structure. 
We use the following approach of Sikorav~\cite{Sikorav}
(or equivalently, Audin-Lalonde-Polterovich~\cite{ALP})
to ensure that the relevant moduli spaces are compact,
and hence the corresponding counts are finite.

First, as explained in \cite{NZ}, the cotangent bundle $T^*X$ equipped
with the canonical symplectic form $\omega$,
conical almost complex structure $J_{con}$, and conical metric $g_{con}$ is tame
in the sense of ~\cite{Sikorav}. To see this, one can verify that $g_{con}$
is conical near infinity, and so it is easy to derive an upper bound
on its curvature and a positive lower bound on its injectivity radius.

Next, given a finite
collection of branes $L_0,\ldots, L_d$, denote by $L$ the union of their perturbations 
$\psi_i( \varphi_{H_i,\delta_i}(L_i))$ as described above.
By construction, the intersection
of $L$ with the region $|\xi|>r$ is a tame submanifold 
(in the sense of ~\cite{Sikorav})
with respect to the structures
$\omega$, $J_{con}$, and $g_{con}$.
Namely, 
there exists $\rho_L>0$ such
that for every $x\in L$, the set of
points $y\in L$ of distance $d_{}(x,y) \leq \rho_L$ is
contractible, and
there exists $C_L$ giving a two-point distance condition
$d_L(x,y) \leq C_L d_{}(x,y)$  whenever $x,y\in L$ with $d_{}(x,y)<\rho_L$.

Now, consider a fixed topological type of
pseudoholomorphic map
$$
u:(D,\del D) \to (T^*X, L).
$$
Assume that all $u(D)$ intersect a fixed compact region,
and there is an a priori area bound ${\rm Area}(u(D))< A$.
Then as proven in \cite{Sikorav}, one has compactness of the moduli space
of such maps $u$. In fact, one has a diameter bound (depending only on the given constants)
constraining how far the image $u(D)$ can stretch from 
the compact set.

In the situation at hand, for a given $A_\infty$-structure constant, we must consider
pseudoholomorphic maps $u$ from polygons with labeled boundary edges.
In particular, all such maps $u$ have image
intersecting the compact set given by a single intersection point. 
The area of the image $u(D)$ can be expressed
as the contour integral 
$$
{\rm Area}(u(D)) = \int_{u(\del D)} \theta.
$$
Since each of the individual Lagrangian branes making up $L$ is exact,
the contour integral only depends upon the integral of $\theta$ along minimal paths between 
intersection points. Thus such maps $u$ satisfy an a priori
area bound. We conclude that
 for each $A_\infty$-structure constant, 
the moduli space defining the structure constant
is compact, and its points are represented by maps $u$ 
with image bounded by a fixed distance from any of the intersection points.

Finally, as usual, the composition map 
$$
m^d:
{hom}_{F(T^*X)}(L_0,L_1)\otimes\dots\otimes{hom}_{F(T^*X)}(L_{d-1},L_d)\rightarrow
{hom}_{F(T^*X)}(L_0,L_d)[2-d]
$$ 
is defined as follows.  Consider elements 
$p_i\in {hom}(L_i,L_{i+1}),$
for $i=0,\ldots,d-1$,
and $p_d\in hom(L_0,L_d)$.
Then the coefficient of $p_d$ in 
$m^d(p_0,\dots,p_{d-1})$
is defined to be the signed sum over pseudoholomorphic maps from a disk
with $d+1$ counterclockwise cyclically ordered
marked points mapping to the $p_i$ and corresponding
boundary arcs mapping to the perturbations of $L_{i+1}.$  Each map contributes
according to the holonomy of its boundary, where adjacent perturbed
components $L_i$ and $L_{i+1}$ are glued with $p_i.$

\medskip

Continuation maps with respect to families of perturbed branes ensure the consistency of
all of our definitions. 
While the details of this were not elaborated on in~\cite{NZ},
Section~\ref{non-char} and the Appendix contain a discussion about continuation maps
which contains what is needed here as a special case.

\medskip

Consider the dg category of right modules over the Fukaya category of $T^*X$.
Throughout this paper, 
we write $F(T^*X)$ for the
the full subcategory
of twisted complexes of representable modules,
and refer to it  as the triangulated envelope of the Fukaya category.
We use the term Lagrangian brane to refer to an object of the Fukaya category,
and brane to refer to an object of its triangulated envelope $F(T^*X)$.

Before continuing, it is worth remarking about the status of units for Lagrangian branes, and thus the precise
relation between the Fukaya category and its triangulated envelope $F(T^*X)$.
Thanks to standard arguments, 
all of the Lagrangian branes considered in this paper 
(in particular, all those arising from sheaves via microlocalization 
as explained in Section~\ref{sect micro}) are cohomologically unital. 
Furthermore, we expect every Lagrangian brane to be cohomologically unital, though
we have not attempted to show this. Rather, our definition of 
the triangulated envelope $F(T^*X)$ as a category of modules
automatically provides (strict) units. The Yoneda embedding from the Fukaya category
to $F(T^*X)$ is cohomologically fully faithful on
cohomologically unital Lagrangian branes. Likewise, since $F(T^*X)$ is unital, the Yoneda embedding
from $F(T^*X)$ to (left or right) modules over $F(T^*X)$ is cohomologically fully faithful. 
This does not rule out the possibility of exotic Lagrangian
branes that are for instance orthogonal to all other branes including themselves.
While potentially interesting, exploring such phenomena is beyond the aims
of this paper.

It is also worth remarking that since $Sh_c(X)$ is split-closed,
as a consequence of our main result, it follows that $F(T^*X)$
is split-closed as well.



\subsection{Duality and time reversal}\label{sec duality}

\subsubsection{Duality}
We introduce here the duality on branes that corresponds to Verdier duality on sheaves.
In Section~\ref{sec duality revisited}, as a consequence of our main result,
we will confirm this compatibility.

Consider the antipodal anti-symplectomorphism
$$
a:T^* X \to T^*X
\qquad a(x,\xi)= (x,-\xi).
$$
It induces a duality equivalence 
$$
\alpha_{X}:F(T^*X)^\circ \risom F(T^*X).
$$
On Lagrangian branes, $\alpha_X$ is given by the map
$$
(L, \CE, \tilde \alpha, \flat) 
\mapsto 
(a(L), a^*(\CE^\vee)\otimes \orient_{X}, -\dim X - a^*(\tilde \alpha), a^*(\flat)).
$$
Here $\orient_{X}$ denotes the pullback to $T^*X$ of the orientation local system of $X$.
Note as well that given a taming perturbation $\psi$ for $L$, one can take the composition
$\psi\circ a$ as a taming perturbation for $a(L)$.
On morphisms and pseudoholomorphic disks, $\alpha_X$ is given by transport
of structure via the antipodal
map $a$.


\subsubsection{Time reversal}
Although our proofs will not require the material of this section, we include
it as a prelude to the informal discussion of Section~\ref{sec corr}.

Let $(T^*X)^-$ denote the cotangent bundle with its opposite symplectic structure.
So except for the symplectic form being negated, no other aspect of the geometry is changed.
In particular, we continue to work with a Riemannian metric on $X$ for which the
notion of normalized geodesic flow is unchanged.

We can repeat the construction of $F(T^*X)$ word for word in order to construct $F((T^*X)^-)$.
So when perturbing branes, we continue to work with Hamiltonian functions $H:(T^*X)^- \to \R$
which are controlled in the sense that there is a real number $r>0$ such that in the region
$|\xi|>r$ we have $H(x,\xi)= |\xi|$. Here the corresponding Hamiltonian isotopy 
$\varphi_{H,t}: (T^*X)^-\to (T^*X)^-$ equals 
normalized anti-geodesic flow $\gamma_{-t}$ in the region $|\xi|>r$
because we are dealing with the opposite symplectic structure.

To calculate the $A_\infty$-structure among
an ordered collection of Lagrangians branes $L_0,\ldots, L_d\subset (T^*X)^-$, 
we continue to repeat our previous definition
and choose controlled Hamiltonian functions $H_0,\ldots, H_d$,
and a fringed set $R\subset \R^{d+1}$
such that for $(\delta_d,\ldots, \delta_0)\in R$, there is a real number $r>0$
such that for any $i\not = j$, we have
$$
\varphi_{H_i,\delta_i}(\ol L_i)\cap 
\varphi_{H_j,\delta_j}(\ol L_j) 
\quad
\mbox{lies in the region $|\xi|<r$.}
$$
By a further compactly supported Hamiltonian perturbation,
we may also arrange so that the intersections are transverse.

Similarly, the rest of the definition of $F((T^*X)^-)$ continues to follow
that of $F(T^*X)$ word by word.

In case of confusion concerning the above perturbation procedure near infinity, 
it is useful to return to the mnemonic that morphisms propagate forward in time.
Thus in the case of $F(T^*X)$, to calculate the morphisms $\hom_{F(T^*X)}(L_0, L_1)$,
we have required that $L_0,L_1$ are perturbed near infinity by normalized geodesic flow
so that $L_1$ is further in the future than $L_0$.
In the case of $F((T^*X)^-)$, we think of the opposite symplectic structure as reversing
the timeline.
To calculate the morphisms $\hom_{F((T^*X)^-)}(L_0, L_1)$,
we again perturb $L_0,L_1$ near infinity 
so that $L_1$ is further in the future than $L_0$. 
But now this implies that near infinity we must perturb $L_1$ 
by normalized anti-geodesic flow
a greater amount than we perturb $L_0$.
Since we use the opposite symplectic structure here,
this implies that we continue to use controlled Hamiltonian functions.

Finally, there is a time reversal equivalence
$$
\rho_X:F(T^*X)^\circ\risom F((T^*X)^-). 
$$
On Lagrangian branes, $\rho_X$ is given by the map
$$
(L, \CE, \tilde \alpha, \flat) \mapsto (L, \CE^\vee, -\tilde\alpha, \flat).
$$
On morphisms and pseudoholomorphic disks, $\rho_X$ is induced by the identity map.


\subsection{Microlocalization}\label{sect micro}

We review here the microlocalization quasi-embedding
$$
\xymatrix{
\mu_X:\Sh_c(X)\ar@{^(->}[r] & F(T^*X)
}
$$ 
constructed in~\cite{NZ}. 
Some useful notation: for a function $m:X\to \R$ and number $r\in R$, we write $X_{m = r}$
for the subset $\{x\in X | m(x) = r\}$ and similarly for inequalities.

\medskip

Let $i: U\hra X$ be an open
submanifold that is a $\CC$-subset of $X$.
Since the complement $X\setminus U$ is a closed $\CC$-subset of $X$, 
we can find a non-negative function
$m:X\to \R_{\geq 0}$ such that $X\setminus U$ is precisely the zero-set of $m$.
Since the complement of the critical values of $m$ form an open $\CC$-subset of $\R$,
 the subset $X_{m >\eta}$
is an open submanifold with smooth hypersurface boundary $X_{m =\eta} $,
for any sufficiently small $\eta >0$.

\medskip

Now let $i_{\alpha}: U_\alpha\hra X$, for $\alpha=0,\ldots, d$,
be a finite collection of open
submanifolds that are $\CC$-subsets of $X$.
Fix non-negative function
$m_\alpha:X\to \R_{\geq 0}$,  for $\alpha=0,\ldots, d$, 
such that $X\setminus U_\alpha$ is precisely the zero-set of $m_\alpha$.
There is a fringed set $R\subset \R^{d+1}_+$ such that 
for any $(\eta_d,\ldots, \eta_0)\in R$, the following holds. First the hypersurfaces 
$X_{m_\alpha =\eta_\alpha}$ are all transverse. Second, for $\alpha<\beta$,
there is a quasi-isomorphism 
of complexes
$$
\hom_{Sh_c(X)}(i_{\alpha*}\C_{U_\alpha}, i_{\beta*}\C_{U_\beta})
\simeq
(\Omega(X_{m_\alpha\geq \eta_\alpha}\cap X_{m_\beta>\eta_\beta}, 
X_{m_\alpha=\eta_\alpha}\cap X_{m_\beta>\eta_\beta}),d)
$$
where $(\Omega, d)$ denotes the relative de Rham complex which calculates the
cohomology of the pair. Furthermore, the composition of morphisms in $Sh_c(X)$
corresponds to the wedge product of forms.

\medskip

Next let 
$f_\alpha:X_{m_\alpha > \eta_\alpha}\to \R$,  for $\alpha=0,\ldots, d$, 
be the logarithm $f_\alpha = \log m_\alpha$.
While choosing the sequence of parameters $(\eta_d, \ldots, \eta_0)$, we can also
choose a sequence of small positive parameters $(\epsilon_d, \ldots, \epsilon_0)$ such that 
the following holds. 
For any $\alpha <\beta$, consider the open submanifold 
$X_{m_\alpha >\eta_\alpha,m_\beta>\eta_\beta} = X_{m_\alpha >\eta_\alpha}\cap 
X_{m_\beta>\eta_\beta} $ with corners equipped with the function
$f_{\alpha, \beta} = \epsilon_\beta f_\beta - \epsilon_\alpha f_\alpha$.
Then there is an open set of Riemannian metrics on $X$ such that for all $\alpha<\beta$,
it makes sense to consider the Morse complex
$\CM(X_{m_\alpha >\eta_\alpha,m_\beta>\eta_\beta}, f_{\alpha,\beta})$,
and
there is a quasi-isomorphism
$$
(\Omega(X_{m_\alpha\geq \eta_\alpha}\cap X_{m_\beta>\eta_\beta}, 
X_{m_\alpha=\eta_\alpha}\cap X_{m_\beta>\eta_\beta}),d)
\simeq
\CM(X_{m_\alpha >\eta_\alpha,m_\beta>\eta_\beta}, f_{\alpha,\beta}).
$$
Furthermore, following arguments of~\cite{HL}, \cite{kontsoib},
homological perturbation theory
provides a quasi-equivalence
between the $A_\infty$-composition structure on the collection of Morse complexes
and the dg structure given by the wedge product of forms.

\medskip

Finally, we 
define the microlocalization quasi-embedding
$$
\xymatrix{
\mu_X:\Sh_c(X)\ar@{^(->}[r] & F(T^*X)
}
$$ 
as follows. Recall by Lemma~\ref{standard star envelope}, that the standard objects
$i_*\C_U$ associated to open submanifolds $i:U\hra X$ generate the constructible
dg derived category $Sh_c(X)$. Thus to construct $\mu_X$, it suffices to find a 
parallel collection of standard objects of $F(T^*X)$.

Given an open submanifold $i:U\hra  X$
and function $m:X\to\R_{\geq 0}$ with zero-set the complement $X\setminus U$,
define
the {standard Lagrangian} $L_{U,f*}\hra T^*X\vert_U $ to be
the graph
$$
L_{U,f*}= \Gamma_{df},
$$
where $df $ denotes the differential of the logarithm $f=\log m$.

The standard Lagrangian $L_{U,f*}$ comes
equipped with a canonical brane structure $(\CE, \tilde\alpha, \flat)$ and taming perturbation $\psi$. 
Its flat vector bundle $\CE$ is trivial, and
its grading $\tilde \alpha$ and relative pin structure $\flat$ are the canonical structures on a graph.
Its taming perturbation $\psi$ is given by the family of standard Lagrangians 
$$
L_{X_{m = \eta}, f_\eta *} = \Gamma_{df_\eta},
\quad
\mbox{ for sufficiently small $\eta>0$},
$$
where $f_\eta =\log m_\eta$ is the logarithm of the shifted function $m_\eta = m -\eta$.

Now one can extend the fundamental result of Fukaya-Oh~\cite{FO} identifying Morse
moduli spaces and Fukaya moduli spaces to the current setting. Namely, one can show
that for any finite ordered collection of open
submanifolds
$i_{\alpha}: U_\alpha\hra X$, for $\alpha=0,\ldots, d$, and any finite collection
of $A_\infty$-compositions respecting the order, there is a fringed set $R\subset \R^{d+1}$
such that for any parameters $(\eta_d,\ldots, \eta_0)\in R$,
the Morse moduli spaces
of the ordered collection of functions $f_{\eta_\alpha}$ are isomorphic
to the Fukaya moduli spaces of the ordered collection of standard branes $L_{X_{m=\eta}, f_\eta*}$
(after further variable dilations of the functions and branes).

Thus we can define a quasi-embedding $\mu_X$ so that on objects we have
$$
\mu_X(i_*\C_U) = L_{U, f*}
$$ 
where $f=\log m$ for any choice of non-negative  
function $m:X\to \R_{\geq 0 }$ such that the complement $X\setminus U$
is the zero-set of $m$. In particular, the standard branes for different choices of $m$
are all isomorphic, and so we will choose one and simply denote it by $L_{U *}$.

\medskip

In what follows, it will also be useful to recall the result of \cite{NZ} describing
where $\mu_X$ takes other standard objects.

Consider the standard sheaf $i_*\C_Y$
associated to an arbitrary submanifold $i:Y\hookrightarrow X$.
Given a non-negative function $m:X\to \R_{\geq 0}$ with zero-set the boundary 
$\del Y =\ol Y\setminus Y$,
define the standard Lagrangian $L_{Y,f*}\hra T^*X|_Y$ to be the fiberwise sum
$$
L_{Y,f*} =T^*_Y X + \Gamma_{df}
$$
where $T_Y^*X\hra T^*X$ denotes the conormal bundle to $Y$, and $\Gamma_{df}\hra T^*X|_Y$
the graph of the differential of the logarithm $f=\log m$.
By construction,
$L_{Y,f*}$ depends only on the restriction
of $m$ to $Y$. 

The standard Lagrangian $L_{Y,f*}$ comes
equipped with a canonical brane structure $(\CE, \tilde\alpha, \flat)$ and taming perturbation $\psi$. 
Its flat vector bundle $\CE$ is the pullback of the normal orientation bundle
$\pi^*(or_X\otimes or_Y^{-1})$,
where $or_X, or_Y$ denote the orientation bundles of $X,Y$ respectively.
Its grading $\tilde \alpha$ is characterized by the following property. Suppose we 
perturb $L_{Y,f*}$ so that it becomes a graph over an open set.
Then if the perturbation is in the direction of anti-geodesic flow near infinity,
the transported grading coincides with the canonical grading carried by the graph.
If the perturbation is in the direction of geodesic flow near infinity,
then the transported grading is equal to the shift by $\codim Y$ of the canonical grading on the graph. Under either such perturbation,
its relative pin structure $\flat$ coincides with the canonical such structure on the graph.
Finally, its taming perturbation $\psi$ is given by the family of standard Lagrangians 
$$
L_{Y_{m = \eta}, f_\eta *} = \Gamma_{df_\eta},
\quad
\mbox{ for sufficiently small $\eta>0$},
$$
where $f_\eta =\log m_\eta$ is the logarithm of the shifted function $m_\eta = m -\eta$.

Then by~\cite{NZ}, the microlocalization 
$\mu_X(i_*\C_Y)$ is isomorphic to the standard brane $L_{Y,f*}$.
In particular, the standard branes $L_{Y,f*}$ for different
choices of $f$ are all isomorphic, and so we will choose one and simply denote it by $L_{Y*}$. 

\medskip

In the next section, it will be helpful to have in mind the following
aspect of the construction of $\mu_X$. 
Our perturbation conventions
allow for all of the standard branes to be perturbed near infinity in the direction of normalized geodesic flow. In particular, this implies that the perturbations contract
the boundaries of standard branes towards the interior of the
corresponding submanifolds.
Thus all calculations can be understood in terms of submanifolds with smooth,
transversely intersecting boundaries.


\subsection{Costandard branes}\label{sect costandard branes}
The material of this section is needed in Section~\ref{sec duality revisited} to confirm
that the microlocalization $\mu_X$ intertwines Verdier duality $\CD_X$ and
the brane duality $\alpha_X$. It does not play a role in the proof in Section~\ref{sect quasi-equivalence} 
that $\mu_X$ is a quasi-equivalence.

\medskip

Let $i:U\hra X$ be an open submanifold, and $i_!\omega_U$ be the corresponding
costandard object.
Let us first study the left Yoneda $\Sh_c(X)$-module
$\CY_\ell(i_{!}\omega_U)$ applied to a finite collection of standard objects
associated to open submanifolds.

For $\alpha=0,\ldots, d$,
consider an open submanifold
$i_{\alpha}: U_\alpha\hra X$, and the corresponding standard object $i_{\alpha*} \C_{U_\alpha}$.
Consider the problem of calculating the directed dg structure among the ordered
collection 
$$i_{!}\omega_U,
 i_{0*}\C_{U_0}, \ldots, i_{d*}\C_{U_d}.
$$

For $\alpha=0,\ldots, d$, fix a
non-negative function
$m_\alpha:X\to \R_{\geq 0}$
such that the complement $X\setminus U_\alpha$ is the zero-set of $m_\alpha$.
Then the techniques of~\cite{NZ}
allow us to fix $i_{!}\omega_U$ but simplify the other objects. To be precise, 
there is a fringed set $R\subset \R^{d+1}_+$ such that 
for any $(\eta_d,\ldots, \eta_0)\in R$, the above directed dg structure 
is quasi-equivalent to that of the ordered collection 
$$
i_{!}\omega_U,
 i_{0,\eta_0*}\C_{X_{m_0>\eta_0}}, \ldots, i_{d,\eta_d*}\C_{X_{m_d>\eta_d}}
$$
where $i_{\alpha,\eta_\alpha} :X_{m_\alpha >\eta_\alpha} \hra X$
denotes the inclusion.

Furthermore, another application 
of the techniques of~\cite{NZ}, along with the adjunction identity
$$
\hom_{Sh_c(X)}(i_{!}\omega_U,i_{\alpha,\eta_\alpha *}\C_{X_{m_\alpha >\eta_\alpha}})
\simeq
\hom_{Sh_c(X)}(\omega_U, i^!i_{\alpha,\eta_\alpha *}\C_{X_{m_\alpha >\eta_\alpha}}),
$$
shows that for sufficently small $\eta>0$,
we can replace the above ordered collection by the ordered collection
$$
i_{\eta !}\omega_{X_{m>\eta}},
 i_{0,\eta_0*}\C_{X_{m_0>\eta_0}}, \ldots, i_{d,\eta_d*}\C_{X_{m_d>\eta_d}},
$$
where as above 
 $i_{\eta} :X_{m >\eta} \hra X$
denotes the inclusion.

Finally, we can calculate the above morphism complexes via 
de Rham complexes. Namely, there are quasi-isomorphisms
of complexes
$$
\hom_{Sh_c(X)}( i_!\omega_{U}, i_{\alpha*}\C_{U_\alpha})
\simeq
(\Omega(X_{m> \eta}\cap X_{m_\alpha>\eta_\alpha}; \omega^{\vee}_X),d),
$$
and for $\alpha<\beta$, there are quasi-isomorphisms
of complexes
$$
\hom_{Sh_c(X)}( i_{\alpha *}\C_{U_\alpha}, i_{\beta*}\C_{U_\beta})
\simeq
(\Omega(X_{m_\alpha\geq \eta_\alpha}\cap X_{m_\beta>\eta_\beta}, 
X_{m_\alpha=\eta_\alpha}\cap X_{m_\beta>\eta_\beta}),d).
$$
Furthermore, the composition of morphisms 
corresponds to the wedge product of forms.

\medskip

Now in parallel with standard branes,
we define costandard branes as follows. For simplicity, and since it suffices
for our later needs, we restrict to the case of open submanifolds.

Given an open submanifold $i:U\hra  X$
and function $m:X\to\R_{\geq 0}$ with zero-set the complement $X\setminus U$,
define
the {costandard Lagrangian} $L_{U,f!}\hra T^*X\vert_U $ to be
the graph
$$
L_{U,f!}= -\Gamma_{df},
$$
where $df $ denotes the differential of the logarithm $f=\log m$.

The costandard Lagrangian $L_{U,f!}$ comes
equipped with a canonical brane structure $(\CE, \tilde\alpha, \flat)$ and taming perturbation $\psi$. 
Its flat vector bundle $\CE$ is the pullback of the orientation bundle $\pi^*(or_X)$,
its grading $\tilde \alpha$ is the shift by $\dim X$ of the canonical grading on a graph, and
its relative pin structure $\flat$ is the canonical structure on a graph.
Finally, its taming perturbation $\psi$ is given by the family of costandard Lagrangians 
$$
L_{X_{m = \eta}, f_\eta !} = -\Gamma_{df_\eta},
\quad
\mbox{ for sufficiently small $\eta>0$},
$$
where $f_\eta =\log m_\eta$ is the logarithm of the shifted function $m_\eta = m -\eta$.

Alternatively, we could use the brane duality $\alpha_X$, and take as definition the
motivating identity
$$
L_{U, f!} \simeq \alpha_X(L_{U, f*}).
$$
In particular, the costandard branes $L_{U,f!}$ for different
choices of $f$ are all isomorphic, and so we will choose one
and denote it by~$L_{U!}$.

\medskip

The construction of the microlocalization $\mu_X$ favors standard objects over costandard objects. 
For example, without further arguments, it does not follow immediately that
$\mu_X$ takes the costandard sheaf $i_{!}\omega_U$ to a costandard brane $L_{U!}$.
Equivalently, without further arguments, it does not immediately follow
that $\mu_X$ intertwines Verdier duality $\CD_X$ with brane duality $\alpha_X$.
To eventually see this (cf. Proposition~\ref{dual confirm}), we will use the following partial identification
of costandard branes.

Consider the microlocalization $\mu_X$, and the resulting pullback functor on left 
modules
$$
\mu^*_X:\lmod(F(T^*X)) \to \lmod(Sh_c(X))
$$
$$
\mu_X^*(\CM)= \CM\circ\mu_X.
$$
In particular, for an object $L$ of $F(T^*X)$, composing $\mu_X^*$ with the Yoneda embedding 
$\CY_\ell$ for left modules provides a left $Sh_c(X)$-module
$$
\mu_X^*(\CY_\ell(L))(\CF) = \hom_{F(T^*X)}(L, \mu_X(\CF)).
$$

\begin{prop}\label{embed costandards}
For any open submanifold $i:U\hra X$, there is a quasi-isomorphism of left $\Sh_c(X)$-modules
$$
\mu_X^*(\CY_{\ell}(L_{U!})) \simeq \CY_\ell(i_{!}\omega_U):Sh_c(X)\to \Ch.
$$
\end{prop}

\begin{proof} 
Let us understand the left $\Sh_c(X)$-module $\mu_X^*(\CY_{\ell}(L_{U!}))$
applied to a finite collection of standard objects
associated to open submanifolds.

Fix a representative $L_{U,f!}$ where as usual $f=\log m$ for a 
non-negative function
$m:X\to \R_{\geq 0}$
such that the complement $X\setminus U$ is the zero-set of $m$.
For $\alpha=0,\ldots, d$,
consider an open submanifold
$i_{\alpha}: U_\alpha\hra X$, a
non-negative function
$m_\alpha:X\to \R_{\geq 0}$
such that the complement $X\setminus U_\alpha$ is the zero-set of $m_\alpha$,
and the corresponding standard brane
$L_{U_\alpha, f_\alpha *}$ where as usual $f_\alpha = \log m_\alpha$,

Consider the problem of calculating the directed $A_\infty$-composition maps among the ordered
collection 
$$L_{U, f!}, L_{U_0, f_0*}, \ldots, L_{U_d, f_d*}.
$$ 
Then our perturbation conventions
allow us to fix $L_{U, f!}$ but perturb the other branes. To be precise, 
there is a fringed set $R\subset \R^{d+1}_+$ such that 
for any $(\eta_d,\ldots, \eta_0)\in R$, we must calculate the 
directed $A_\infty$-composition maps among the ordered
collection 
$$L_{U, f!}, L_{X_{m_0 >\eta_0}, f_{0, \eta}*}, \ldots, L_{X_{m_d>\eta_d}, f_{d,\eta}*}.
$$
Here as earlier the underlying Lagrangian of
the standard brane $L_{X_{m_\alpha > \eta_\alpha}, f_{\alpha, \eta}*}$
is the graph of the differential of 
$f_{\alpha, \eta_\alpha} = \log (m_\alpha  -\eta_\alpha)$
over the open set $X_{m_\alpha >\eta_\alpha}$ 

By definition, the taming perturbation of the costandard brane $L_{U, f!}$ moves it to the
costandard brane $L_{X_{m>\eta}, f_\eta !}$, for sufficiently small $\eta>0$.
Here as above the underlying Lagrangian of the costandard brane $L_{X_{m>\eta}, f_\eta !}$
is the negative of the graph of the differential of 
$f_{\eta} = \log (m  -\eta)$
over the open set $X_{m >\eta}$.
Thus we are left to calculate the 
directed $A_\infty$-composition maps among the ordered
collection 
$$L_{X_{m>\eta}, f_\eta !}, 
L_{X_{m_0 >\eta_0}, f_{0, \eta}*}, \ldots, L_{X_{m_d>\eta_d}, f_{d,\eta}*}.
$$

The techniques of~\cite{NZ} extend directly to this situation: the relevant
Fukaya moduli spaces can be identified with the corresponding Morse moduli spaces.
In turn, following arguments of~\cite{HL}, \cite{kontsoib},
homological perturbation theory
provides a quasi-equivalence
between the Morse $A_\infty$-composition structure and the dg structure given by 
the wedge product of forms. Finally, as discussed above, the dg structure on differential
forms calculates the dg structure on the corresponding constructible sheaves.
\end{proof}


\subsection{Non-characteristic isotopies}\label{non-char}

We discuss here the invariance of calculations among microlocal branes under 
a very specific class of non-characteristic 
Hamiltonian isotopies. What we explain is the minimum technical result needed
to establish our main theorem. Further generalizations are discussed in the Appendix.

\subsubsection{Motivation: sheaf calculations}
This section is intended as motivation for the Floer calculations to follow, but
is not logically needed in what follows.

Let's consider a family of stratifications
of $X$ parametrized by the real line $\R$.
More precisesly, by a one-parameter family of stratifications of $X$, 
we mean a single Whitney stratification
$
\mathfrak S =\{\fS_\beta\}
$
of  
$
\R\times X
$
satisfying the following:

\begin{enumerate}

\item 
The restrictions of the projection $p_\R: \R\times X\to \R$
to each stratum $\fS_\beta$ of $\fS$ is nonsingular.

 \item There is a compact interval $[a,b]\hra \R$ such that 
 the induced stratification of 
 $p_\R^{-1}(\R\setminus [a,b])$ obtained by restricting $\fS$
 is locally constant.

\end{enumerate}

For each $s\in \R$, we will denote by $\fS(s)= p_\R^{-1}(s) \cap \fS$ the fiber of $\fS$.
Note that condition (1) above implies that the topological type of the stratification $\fS(s)$
is constant with respect
to $s\in\R$.
More precisely, by the Thom Isotopy Lemma,
one can construct a homeomorphism $\psi:\R\times X\to \R\times X$
such that $p_\R\circ\psi = p_\R$, $\psi(\fS) = \R\times \fS(0)$.

\medskip

Suppose two one-parameter family of stratifications $\fS =\{\fS_\beta\}$, $\fS'=\{\fS'_\alpha\}$ 
of $X$ are transverse. Note that this is equivalent to the fibers $\fS(s)$ and $\fS'(s)$
being transverse for all $s\in\R$.
Then the
Whitney stratification $\fS\cap \fS'$
of $\R\times X$
with strata the intersections of strata $\{\fS_\beta \cap \fS'_\alpha\}$
is again a one-parameter family of stratifications of $X$.

In particular, 
fix a Whitney stratification $\CS=\{S_\alpha\}$ of $X$, and let $\CS_\R=\{\R\times S_\alpha\}$
be the constant one-parameter family of stratifications of $\R\times X$.
We will say that $\fS$ is $\CS$-non-characteristic if $\fS$ is
transverse to $\CS_\R$. 

\medskip

Now consider the $\CS$-constructible dg derived category $Sh_\CS(X)$.
Consider as well any object $\CF$ of the $\fS$-constructible dg derived category $Sh_{\fS}(\R\times X)$,
and denote by $\CF_s$ its restriction to the fiber $X=p_\R^{-1}(s)$.
In general, $\CF_s$ is not an object of $Sh_{\CS}(X)$, 
but can be paired with objects of $Sh_\CS(X)$
in the ambient constructible dg derived category $Sh_c(X)$.

\begin{lem}\label{lem sheaf non-char}
Suppose $\fS$ 
is an $\CS$-non-characteristic one-parameter family of stratifications
of $X$.
Then for any object $\CF$ of $Sh_{\fS}(\R\times X)$,
and any test object $\CP$ of $Sh_{\CS}(X)$, there are functorial quasi-isomorphisms among
the complexes
$$
\hom_{Sh_c(X)}(\CP, \CF_s),
\qquad
\mbox{for all $s\in\R$}.
$$
\end{lem}

\begin{proof}
Consider the object of $Sh_c(\R)$ given by the pushforward
$$
p_{\R!}\intHom_{Sh_c(\R\times X)}(\C_\R\boxtimes \CP, \CF).
$$
By base change, its stalk at $s\in\R$ is the complex $\hom_{Sh_c(X)}(\CP, \CF)$.
As mentioned above, the Thom Isotopy Lemma provides a stratum-preserving
homeomorphism
$$\psi:\R\times X\to \R\times X
\qquad
\psi( \fS \cap \CS_\R) = \R\times (\fS(0) \cap \CS)
$$ 
such that $p_{\R}\circ\psi= p_\R$.
Thus the projection via $p_{\R!}\simeq p_{\R!}\circ\psi_!$
of any object of $Sh_{\fS\cap \CS_\R}(\R\times X)$ is constant.
\end{proof}

Our aim in the next two sections is to produce an analogue of Lemma~\ref{lem sheaf non-char}
with sheaves on $X$ replaced by branes in $T^*X$.
In particular, we will be interested in a brane version of the following specific example.

\begin{ex} Let $\fS$ be
a one-parameter family of stratifications of $X$. Suppose further that $\fS$ consists of three strata:
a smooth submanifold
$$
\fraki:\fY \hra \R\times X,
$$
its smooth boundary $\del\fY =\ol \fY\setminus \fY\hra \R\times X$, and their complement $(\R\times X)\setminus \ol \fY$.
For each $s\in \R$, we will denote by 
$$
\fraki_s:\fY_s = p_\R^{-1}(s) \cap \fY \hra X
$$
 the fiber of $\fY$.

For a given stratification $\CS$ of $X$,
we will say that $\fY$ is an $\CS$-non-characteristic one-parameter family
of submanifolds with smooth boundaries if the corresponding three stratum stratification $\fS$ 
is $\CS$-non-characteristic.

If $\fY$ is an $\CS$-non-characteristic
one-parameter family
of submanifolds  with smooth boundaries,
then 
for any test object $\CP$ of $Sh_{\CS}(X)$, there are functorial quasi-isomorphisms among
the complexes
$$
\hom_{Sh_c(X)}(\CP, \fraki_{s*}\C_{\fY_s}),
\qquad
\mbox{for all $s\in\R$}.
$$
\end{ex}


\subsubsection{Continuation of Floer calculations}

Let $L\hra T^*X$ be a tame exact Lagrangian brane. Consider a
time-dependent Hamiltonian function $H_s: X\times \R\to \R$
such that its differential $dH_s$ is  {\em compactly supported} in $X$ and $\R$.
Let $\varphi_s:T^*X\to T^*X$ be the associated
Hamiltonian flow. 
Acting on the initial brane $L$,
we obtain a family of tame exact Lagrangian
branes $\fL_s = \varphi_s(L)$
satisfying the following:

\begin{enumerate}

\item $\fL_0 = L$.

\item  Near infinity in $\R$, the family $\fL_s$ is locally constant: there is a compact interval $[a,b]\hra \R$ such that for $s\in \R\setminus [a,b]$, 
the family $\fL_s$ is locally constant.

\item Near infinity in $T^*X$, the family $\fL_s$ is constant: there exists $r>0$ such that
$$
\{|\xi|>r\}\cap \fL_s = \{|\xi|>r\}\cap \fL_0, \mbox{ for all $s\in \R$.}
$$

\end{enumerate}
The family $\fL_s$ is the most general possible compactly supported motion
of the exact brane $L$.
We will consider more general families in the next section.

The following is a straightforward generalization of by now standard techniques
in Floer theory. Namely, once we confirm the necessary a priori estimates on the 
possible diameters of
pseudoholomorphic disks, continuation maps provide the sought after natural
transformations.

\begin{prop}\label{prop compact isotopy}
For any test  $P$ of $F(T^*X)$, there are functorial quasi-isomorphisms among the 
Floer complexes
$$
\hom_{F(T^*X)}(P, \fL_{s}),
\qquad
\mbox{ for all $s\in\R$}.
$$
\end{prop}

\begin{proof}
Fix parameters $a\ll 0\ll b\in \R$.
Following Seidel~\cite[Section 10c]{Seidel},
given test objects $P_1,\ldots, P_d$ of $F(T^*X)$,
we would like to define maps
$$
T_{d}:\hom_{F(T^*X)}(P_1, \fL_a) \otimes 
\left(
\bigotimes_{k=1}^{d-1}
\hom_{F(T^*X)}(P_{k+1},P_k)
\right)
\to 
\hom_{F(T^*X)}(P_d, \fL_b)[1-d]
$$
assembling into an $A_\infty$-transformation of right Yoneda modules. 
For compact branes, one immediately has
the sought after maps: a signed count of pseudoholomorphic disks
with moving boundary conditions given by the family $\fL_s$, for $s\in \R$, and static conditions
given by the test branes $P_1,\ldots, P_d$ provides the structure
constants of the maps. Furthermore, the transformations satisfy a compatibility
with respect to the concatenation of families. In particular, 
the transformations are quasi-isomorphisms since the constant family gives 
the identity functor.

\medskip

In our current setting, to implement this approach, 
we need to be careful to make sure
the relevant moduli spaces remain compact.
This is the content of the rest of the proof of the proposition.

First, following \cite{Seidel}, it is convenient to recast the moduli problem of pseudoholomorphic polygons
with moving boundary conditions in terms of pseudoholomorphic sections
of varying almost complex targets over holomorphic polygons.

Given a non-negative integer $d$,
consider the manifold with boundary $D_d\hra \C$ obtained from the closed unit disk
$\ol D\hra \C$
by removing the $(d+1)$st roots of unity from the boundary $\del D\hra \C$.
Thus the boundary $\del D_d\hra D_d$ 
is a disjoint union of $d+1$ cyclically ordered open intervals $I_k$ 
indexed by $k = 0,1,\ldots, d$. 
We equip the interior $D^\circ_d = D_d\setminus \del D_d$ 
with the standard symplectic structure $\omega_d$ obtained by restricting
the standard exact symplectic structure of $\C$.
(The interaction 
of $\omega_d$ with the boundary $\del D_d$ will
not play a significant role.)
We will consider $D_d$ together with $\omega_d$ as a fixed symplectic manifold with boundary.

Consider the product manifold
$D_d\times T^*X$ with the obvious projection onto the first factor
$
\pi:D_d\times T^*X\to D_d.
$
Identify the boundary 
component $I_0\hra \del D_d$ with the real line $\R$, and
consider the embedded submanifold
$$
\fL_{move} = I_0  \times_{D_d} \fL_s \hra D_d \times T^*X.
$$
Similarly, given a collection of static test objects $P_1,\ldots, P_d\hra T^*X$,
consider the embedded submanifolds
$$
 P_{k,stat} = I_k  \times P_k \hra D_d \times T^*X,
\quad\mbox{ for $k=1,\ldots, d$}.
$$

Consider a one-form $\kappa$ on the base $D_d$ with values in functions
on the fiber $T^*X$. Given a vector field $v$ on $D_d$, the evaluation $\kappa(v)$
provides a family of functions on $T^*X$ parametrized by $D_d$. By passing
to the corresponding Hamiltonian vector fields, we may think of $\kappa$ as a section of the vector bundle $\Hom(TD_d,\pi_*TT^*X)$. 
In other words, $\kappa$ provides a connection $\nabla_\kappa$ on the trivial
family of symplectic manifolds defined by $\pi$.

Now, given the submanifold $\fL_{move}\hra D_d\times T^*X$, 
we will choose a compactly supported one-form $\kappa$ so
that $\fL_{move}$ 
is preserved by the parallel transport of the connection $\nabla_\kappa$. 
To this end, it is convenient to choose a collared neighborhood $\CN_0\hra D_d$
of the boundary component $I_0\hra \del D_d$. 
Thus we have an identification $\CN_0\simeq I_0 \times [0,1)$
which we think of as being given by coordinates $(s,t)$. 
Then it is straightforward to construct a one-form $\kappa$ as above satisfying the following: 
\begin{enumerate}
\item
There is a compact subset of $\CN_0\times T^*X$ outside of which $\kappa$ vanishes.
\item
In local coordinates, $\kappa$ can be expressed as $f(s,t,x,\xi) ds$, where $(x,\xi)$ are local
coordinates on $T^*X$.
\item
Along the boundary component $I_0$, the submanifold $\fL_{move}$ 
is preserved by the parallel transport of the connection $\nabla_\kappa$.
\end{enumerate}

Next, let $j$ be a compatible complex structure on the symplectic manifold $D_d$ (for a surface, any correctly oriented complex structure is compatible), and let $J_{con}$ be the conical almost complex structure
on $T^*X$ discussed in the previous section. Together with $\kappa$, the two structures provide an
almost complex structure $J_\kappa$ on the family $D_d\times T^*X$ such that the projection $\pi$
is pseudoholomorphic. Namely, one takes
$$
J_\kappa(v, w) = (j(v), J_{con}(w) + J_{con}(\sigma_\kappa(j(v))) + \sigma_\kappa(v))
$$
where $\sigma_\kappa\in \Hom(TD_d,\pi_*TT^*X)$ is the section associated to $\kappa$. 
In what follows, we will always consider $D_d\times T^*X$
equipped with the almost complex structure $J_\kappa$. Note that while we will think
of $\kappa$ and $J_{con}$ as fixed, we will allow the choice
of $j$ to vary.

The moduli problem of $J_{con}$-holomorphic polygons in $T^*X$
with moving boundary condition $\fL_s$ and static boundary conditions $P_1,\ldots, P_d$,
coincides
with that of $J_\kappa$-holomorphic sections of $\pi$ with boundary conditions
$\fL_{move},P_{1,stat},\ldots, P_{d,stat}$.
To get bounds for solutions, we introduce a symplectic structure on $D_d\times T^*X$ as follows.
Let $\theta$ be the canonical one-form on $T^*X$, and consider the two-form 
$$
\omega_\kappa = d\theta + d\kappa + c\omega_d
$$
where $c>0$ is some fixed constant.
Using the explicit form of $\kappa$, it it simple to check that we can choose $c$ large enough so
that $\omega_\kappa$ will be non-degenerate.
Furthermore, it is simple to check that the complex structure $J_\kappa$ is compatible with $\omega_\kappa$,
and the boundary conditions
$\fL_{move},P_{1,stat},\ldots, P_{d,stat}$ are Lagrangian.

\medskip

With the preceding setup in hand, 
we are now in a context where we can appeal to standard results to verify the proposition.
By construction, all of the branes under consideration are tame, and so we have a priori
bounds on the diameters of the relevant pesudoholomorphic disks.
Thus standard techniques as outlined in \cite[Section 10c]{Seidel} provide
continuation maps giving a functorial quasi-isomorphism.
\end{proof}


\subsubsection{Non-characteristic isotopies of branes}
In this section, we will consider more general families of exact branes in $T^*X$ parametrized by the real line $\R$.

By a
one-parameter family of closed (but not necessarily compact) 
submanifolds (without boundary) in $T^*X$, we mean a closed 
 submanifold
$$
\fL \hra \R\times T^*X
$$ 
satisfying the following:

\begin{enumerate}

\item 
The restriction of the projection $p_\R: \R\times X\to \R$
to the submanfold $\fL$ is nonsingular.

\item There is a real number $r>0$, such that
the restriction of the product
$
p_\R\times|\xi| : T^*X \to \R \times (r,\infty)
$
to the subset $\{|\xi|>r\}\cap \fL$ is proper and nonsingular. 

 \item There is a compact interval $[a,b]\hra \R$ such that
 the restriction of the projection $p_X:\R\times T^*X\to T^*X$ 
 to the submanifold $p_\R^{-1}([\R\setminus [a,b])\cap \fL$
 is locally constant.

\end{enumerate}

Note that conditions (1) and (2) will be satisfied if 
the restriction of the projection $\ol p_\R:\R\times \ol T^*X \to \R$ 
to the closure $\ol\fL\hra \ol T^*X$ is 
nonsingular as a stratified map, but the weaker condition stated is a useful generalization.
It implies in particular that the fibers $\fL_s = p_\R^{-1}(s) \cap\fL \hra T^*X$ are all diffeomorphic, but imposes no requirement
that their boundaries at infinity should all be homeomorphic as well.

\medskip

By a one-parameter family of tame Lagrangian branes in $T^*X$, we mean
a one-parameter family of closed submanifolds $\fL\hra \R\times X$ in the above sense
such that the fibers $\fL_s = p_\R^{-1}(s) \cap\fL \hra T^*X$ also satisfy:

\begin{enumerate}

\item The fibers $\fL_s$ are exact tame Lagrangians with respect to the usual
symplectic structure and any almost complex structure conical near infinity.

\item The fibers $\fL_s$ are
equipped with a locally constant brane structure $(\CE_s, \tilde\alpha_s,\flat_s)$
with respect to the usual background classes. 

\end{enumerate}

Note that if we assume that $\fL_0$ is an exact Lagrangian, then $\fL_s$ being an
exact Lagrangian
is equivalent to the family $\fL$ being given by the flow $\varphi_{H_s}$
of the vector field of a time-dependent
Hamiltonian $H_s:T^*X\to \R$.
Note as well that a brane structure consists of topological data, so can be transported unambiguously
along the fibers of such a family. 

\medskip

\begin{rmk}
In fact, it is possible to prove continuation of Floer homology for even more general
families.
A motivation for this level of generality is that it allows one to check that all taming perturbations 
for a brane lead to equivalent calculations. 
See the Appendix for a discussion in this direction.
\end{rmk}

\medskip

Fix a conical Lagrangian $\Lambda\subset T^*X$, and let $\Lambda^\infty=\ol \Lambda\cap T^\infty X$ 
be its boundary at infinity. Let $F(T^*X)_\Lambda$ be the full
subcategory of $F(T^*X)$ generated by Lagrangian branes $L$ 
whose boundary at infinity $L^\infty=\ol L\cap T^\infty X$ lies in $\Lambda^\infty$. 

Suppose $\fL\hra \R\times T^*X$ is a one-parameter family of tame Lagrangian branes.
We will say that $\fL$ 
is $\Lambda$-non-characteristic if 
$$
\ol \fL_s\cap \Lambda^\infty =\emptyset,
\qquad
\mbox{ for all $s\in\R$.}
$$

\begin{prop}\label{prop brane non-char}
Suppose $\fL\hra \R\times T^*X$ is a $\Lambda$-non-characteristic
one-parameter family of tame Lagrangian branes.
For any test object $P$ of $F_\Lambda(T^*X)$, there are functorial quasi-isomorphisms among the 
Floer complexes
$$
\hom_{F(T^*X)}(P, \fL_{s}),
\qquad
\mbox{ for all $s\in\R$}.
$$
\end{prop}

\begin{proof}

Our strategy will be to 
``factor" the family $\fL$ into many small steps which each fall into a 
broad class of manageable moving boundary conditions. 
Over each step, we will be able to establish quasi-isomorphisms of right Yoneda modules. 
To prove the proposition, we will then take compositions
of these quasi-isomorphisms.

To begin, since we assume that $\fL$ 
is locally constant away from a compact interval $[a,b]\hra \R$,
it suffices to show that for each fixed
parameter $s_0\in \R$, there is a
neighborhood $I_{s_0}=[a_0,b_0]\subset \R$ with $a_0< s_0< b_0$,
such that there are functorial quasi-isomorphisms among the 
Floer complexes
$$
\hom_{F(T^*X)}(P, \fL_{s}),
\qquad
\mbox{ for all $s\in I_{s_0}$}.
$$

Fix a parameter $s_0\in\R$.
Consider a finite number of $A_\infty$-operations
among $\fL_{s_0}$ and a fixed collection of test objects $P_1,\ldots, P_d$ of $F_\Lambda(T^*X).$
Consider the moduli problem of pseudoholomorphic disks that calculate the structure
constants of the operations.
Recall that we have an a priori diameter bound
on such pseudoholomorphic disks: there is a large constant $r_{s_0}>0$
such that none of the relevant disks enter the region of $T^*X$ given by $|\xi|>r_{s_0}/2$.
Moreover, by the continuity of the diameter bound, we can find an open interval $K_{s_0}\subset \R$
containing $s_0$ such that the pseudoholomorphic disks that calculate
the same structure constants for
$
 \fL_{s},
$
 for all $s\in K_{s_0}$,
do not enter the region of $T^*X$ given by $|\xi|>r_{s_0}$.

\medskip

Now for fixed constants $r_2>r_1> r_0 >r_{s_0}$, 
there exists a sufficiently small closed subinterval $I_{s_0}=[a_{0}, b_0]\subset K_{s_0}$,
with $a_0 < s_0<b_0$,
such that we can ``factor" the family $\fL$ over the parameters $I_{s_0}$ into two families of the following form.

First, we can define a family $\fL'$ over $I_{s_0}$ satisfying the following.
Let $L'\hra I_{s_0} \times T^*X$ be the union of the moving brane 
$\fL'$ and the static branes $I_{s_0}\times P_1,\ldots, I_{s_0}\times P_d$. 
Then by construction, we can arrange that
\begin{enumerate}
\item $\fL'_{a_0}= \fL_{a_0}$;
\item $\fL'$ is constant in the region $|\xi|<r_1$;
\item $\fL'$  coincides with $\fL$ in the region $|\xi|>r_2$;
\item in the region $|\xi| >r_0$, the product function
$$
p\times|\xi| : L' \to \R \times (r_0,\infty)
$$
is nonsingular.
\end{enumerate}

Second, we can define a family $\fL''$ over $I_{s_0}$ satisfying the following
\begin{enumerate}
\item $\fL''_{a_0} = \fL'_{b_0}$;
\item $\fL''$ is constant in the region $|\xi|>r_2$;
\item $\fL''$  coincides with $\fL$ in the region $|\xi|<r_1$;
\item $\fL''_{b_0}=\fL_{b_0}$.
\end{enumerate}

It is worth commenting that our choice to work locally near a fixed parameter $s_0$ is due to the fact that it is easy to find such a factorization locally in $s$.

\medskip

Now if we can establish individually that the Yoneda modules associated to the fibers of such 
families $\fL',\fL''$
are quasi-isomorphic, then it will
follow by composition that the Yoneda modules associated to the fibers of 
the family $\fL$ itself
over the interval $I_{s_0}$ are quasi-isomorphic.

\medskip

\noindent{\em Case 1 ($\fL'$ constant away from infinity).} 
The first case is particularly easy in that there is in fact
a strict identification of the $A_\infty$-operations under consideration.
By construction, for the parameter $s_0\in\R$, we have arranged so that
the pseudoholomorphic disks that calculate the corresponding structure constants
do not leave the the region $|\xi|<r_1$.
Furthermore, we have arranged so that for any $s\in I_{s_0}$, the intersection of $\fL_s$ and the region
$|\xi|<r_1$ is constant, and in particular, the diameter bound constraining the relevant disks
holds independently of $s$.
Thus the identity map on intersection points gives an identification
of the $A_\infty$-operations. 

\medskip
\noindent{\em Case 2 ($\fL''$ constant near infinity).}
With the assumptions of the second case, we must check
that a compactly supported family leads to a quasi-isomorphism of Yoneda modules.
This was the content of Proposition~\ref{prop compact isotopy}.
\end{proof}

A further generalization of Proposition~\ref{prop brane non-char} will be presented
in the Appendix. For now, we will simply mention the single application of 
Proposition~\ref{prop brane non-char} that will be used in what follows.

\begin{ex}\label{ex brane non-char}

Fix a Whitney stratification $\CS=\{S_\alpha\}$ of $X$, 
and let $\Lambda_\CS = \cup_\alpha T^*_{S_\alpha} X\hra T^*X$ be the associated
conical Lagrangian.
Let $\fY\hra \R\times X$ be an $\CS$-non-characteristic one-parameter family of submanifolds with smooth boundaries in $X$.

Fix a non-negative function $m:\R\times X\to\R$
that vanishes precisely on the boundary $\del \fY\hra \R\times X$, and consider the function
$f:\fY\to \R$ given by $f=\log m$. For each $s\in\R$, consider the function $f_s:\fY_s\to \R$ 
obtained by restricting
$f$ to the fiber $\fY_s = p_\R^{-1}(s)\hra X$.

Define
the one-parameter family of closed submanifolds $L_{\fY,f*}\hra \R\times T^*X$ to be
the union of the fiberwise sums
$$
\fL_{\fY,f*}= \bigsqcup_{s\in\R} (s, T^*_{\fY_s} X + \Gamma_{df_s}) \hra\R\times T^*X
$$ 
where $T^*_{\fY_s}  X$ denotes the conormal bundle to $\fY_s\hra  X$,
and $\Gamma_{df_s}$ the graph
of the differential of $f_s$ over $\fY_s$.

By construction, $\fL_{\fY,f*}$ is a $\Lambda_\CS$-non-characteristic one-parameter family
of tame Lagrangian branes.
Thus by Proposition~\ref{prop brane non-char}, for any object $P$ of $F_{\Lambda_\CS}(T^*X)$, there are functorial quasi-isomorphisms among the 
Floer complexes
$$
\hom_{F(T^*X)}(P, \fL_{\fY_s,f_s*}),
\qquad
\mbox{ for all $s\in\R$}.
$$

\end{ex}


\section{Microlocalization is a quasi-equivalence}\label{sect quasi-equivalence}
In this section, we prove that the constructible dg derived category $Sh_c(X)$
is quasi-equivalent to the derived Fukaya category $F(T^*X)$. We show that every object of $F(T^*X)$
is quasi-isomorphic to the microlocalization of an object of $Sh_c(X)$. 



\subsection{Statement of results}
Consider the microlocalization quasi-embedding
$$
\xymatrix{
\mu_X: Sh_c(X) \ar@{^(->}[r] & F(T^*X).
}
$$
By definition, the fact that it is a quasi-embedding means that
it induces a fully faithful functor on cohomology categories
$$
\xymatrix{
H(\mu_X): D_c(X) \ar@{^(->}[r] & DF(T^*X).
}
$$
Thus to show that $\mu_X$ is a quasi-equivalence, we must show
that $H(\mu_X)$ is essentially surjective, or in other words, 
that every object of $DF(T^*X)$
is isomorphic to an object coming from $D_c(X)$.

\medskip

By construction, the image of $\mu_X$ is generated by the standard branes $L_{Y*}\hra T^*X$
associated to submanifolds $Y\hra X$. 
Recall the Yoneda embedding into right modules
$$
\CY_r:F(T^*X)\to \rmod(F(T^*X))
\qquad
\CY_r(L): L' \mapsto \hom_{F(T^*X)}(L', L).
$$
The main technical result of this section is the folllowing.

\begin{thm}\label{theorem projection}
Let $L$ be an object of $F(T^*X)$. The Yoneda module $\CY_r(L)$
can be expressed as a twisted complex of the Yoneda 
modules of standard branes $\CY_r(L_{Y*})$.
\end{thm}

\begin{rmk}
Along the way, we will also directly establish the (weaker) statement that if
$\hom_{F(T^*X)}(L_{\{x\}}, L)$ is acyclic for the standard
branes $L_{\{x\}}\hra T^*X$ associated to all points $x\in X$, then $\hom_{F(T^*X)}(P, L)$
is acyclic for all test objects $P$. 

Complexes of the form $\hom_{F(T^*X)}(L_{\{x\}}, L)$
will appear as the coefficients of the modules $\CY_r(L_{Y*})$ in the decomposition
of the module $\CY_r(L)$. See Remark~\ref{rem coeffs} for a more precise statement
on the structure of the coefficients.
\end{rmk}

Theorem~\ref{theorem projection} immediately implies the following.

\begin{thm}\label{theorem quasi-equivalence}
Microlocalization is a quasi-equivalence
$$
\xymatrix{
\mu_X: Sh_c(X) \ar[r]^-{\sim} & F(T^*X)
}$$
\end{thm}

%

%

The next four sections are devoted to the proof of Theorem~\ref{theorem projection}.
In the final section, we discuss constructibility properties of the quasi-equivalence.


\subsection{Two Floer calculations}

Let $X_0, X_1$ be compact real analytic manifolds.

\medskip

We consider here the product manifold $X_0\times X_1$, and the Fukaya category
of its cotangent bundle $T^*(X_0\times X_1)$ with respect to the usual background structures
(cf. Section~\ref{sec branes} or \cite{NZ}).
Our aim is to identify the Yoneda modules associated to some
simple but important examples of branes on the product.


\subsubsection{Product branes}

Consider the product map on the set of branes
$$
Ob (F(T^*X_0))\times Ob(F(T^*X_1)) \to Ob(F(T^*X_0\times T^*X_1))
$$
$$
((L_0,\CE_0,\tilde\alpha_0,\flat_0),
(L_1,\CE_1,\tilde\alpha_1,\flat_1))
\mapsto 
(L_0\times L_1,\CE_0\otimes \CE_1,\tilde\alpha_0+\tilde\alpha_1,\flat_{0,1}).
$$
Here $\flat_{0,1}$ denotes the canonical induced relative pin structure on $L_0\times L_1$
with respect to the usual background class
$$
(\pi_0 \times \pi_1)^*w_2(X_0\times X_1)  = \pi_0^*w_2(X_0) + 
\pi_0^*w_1(X_0)\cdot\pi_1^*w_1(X_1) + \pi_1^*w_2(X_1).
$$
To construct $\flat_{0,1}$,
observe that the relative pin structures $\flat_0,\flat_1$ with respect to the usual background
classes $\pi_0^*w_2(X_0), \pi_1^*w_2(X_1)$ respectively together provide
a twisted lift of $TL_0\times TL_1$ to the pin group with respect to the background class 
$$
\pi_0^*w_2(X_0) + w_1(L_0)\cdot w_1(L_1) + \pi_1^*w_2(X_1).
$$
By assumption, the Maslov classes of $L_0,L_1$ vanish, and so $w_1(L_0),w_1(L_1)$
are the restrictions of $\pi_0^*w_1(X_0), \pi_1^*w_1(X_1)$ respectively.

We use the term product branes to refer to objects of $F(T^*X_0\times T^*X_1)$ that arise via the preceding construction.
When there is no chance of confusion, we will denote a product brane
by the product of the underlying Lagrangians.

\begin{lem}\label{lem product floer calc}
For any test objects $L_0, P_0$ of $F(T^*X_0)$, and $L_1,P_1$ of $F(T^*X_1)$,
there is a functorial (in each object) quasi-isomorphism of Floer complexes
$$
\hom_{F(T^*X_0\times T^*X_1)}(L_0\times L_1, P_0\times P_1)
\simeq
\hom_{F(T^*X_0)}(L_0, P_0)
\otimes  
\hom_{F(T^*X_1)}(L_1, P_1).
$$
\end{lem}

\begin{proof}
One can choose all further necessary structures such as perturbations to be
a product of the corresponding structures on the factors. In this way, one obtains in fact 
a strict isomorphism of complexes.
\end{proof}


\subsubsection{Diagonal brane}
Consider the smooth, closed submanifold
given by the diagonal $\Delta_X \subset X\times X$,
and let $\C_{\Delta_X}$ denote the constant sheaf along $\Delta_X$.

Let $L_{\Delta_X}$ be the standard object of $F(T^*X \times T^*X)$
obtained as the microlocalization
$$
L_{\Delta_X} \simeq \mu_{X\times X}(\C_{\Delta_X}).
$$
We will refer to $L_{\Delta_X}$ as the diagonal brane though
its underlying Lagrangian is the conormal bundle
$$
T^*_{\Delta_X} (X\times X) =\{(x, \xi; x,-\xi)\in T^*X\times T^*X\}
$$
(so strictly speaking, not truly the diagonal -- the true diagonal is not Lagrangian).
Its flat vector bundle is the pullback of the normal orientation bundle
$or_{X\times X} \otimes or^{-1}_{\Delta_X}$.

\begin{prop}\label{diag}
For any test objects $P_0,P_1$ of $F(T^*X)$, there is a 
functorial quasi-isomorphism of Floer complexes
$$
\hom_{F(T^*X)}(P_1,P_0)\simeq \hom_{F(T^*X\times T^*X)} (L_{\Delta_X}, P_0\times\alpha_X(P_1)).
$$
\end{prop}

\begin{proof}
Consider the intersection points and pseudoholomorphic disks 
involved in calculating the right hand side
as a function of the branes $P_0, P_1$.
Observe that to make any calculation involved, we can fix the brane $L_{\Delta_X}$
and work with perturbations that only move the other branes.

Now apply the product $\id\times a_1$ 
of the identity map $\id$ of the first factor
and the antipodal map of the second factor
$$
a_1:T^* X \risom T^*X 
\qquad 
a_1(x_1,\xi_1) = (x_1, -\xi_1) 
$$
to the objects under consideration.
Observe that
$\id\times a_1$  takes the conormal Lagrangian $T^*_{\Delta_X} (X\times X)$ 
to the diagonal submanifold
$\Delta_{T^*X}$, the Lagrangian $P_0$ to itself,
and the dualized Lagrangian $\alpha_X(P_1)$ back to $P_1$.

Standard gluing arguments imply that $\id\times a_1$  takes a pseudoholomorphic map 
$$
(u_0, u_1):D\to T^*X\times T^*X
$$
with a $L_{\Delta_X}$-labelled boundary component 
$$
(u_0, u_1)|_C:C\to L_{\Delta_X}
$$
to a pseudoholomorphic map 
$$
u_0 \cup a_1(u_1):D\cup_C \ol D \to T^*X
$$
where $\ol D$ denotes the conjugate disk, and $D\cup_C \ol D$ the gluing of $D$, $\ol D$ along $C$.

Tracing through brane structures, we see that the above identification of moduli spaces 
provides the sought-after functorial quasi-isomorphism
$$
\hom_{F(T^*X)}(L_{\Delta_X}, P_0\times\alpha_X(P_1))\simeq
\hom_{F(T^*X)}(P_1,P_0)
$$ 
Note that the appearance of the orientation bundle $or_X$
on the dualized brane $\alpha_X(P_1)$
matches up with the appearance  
of the normal orientation bundle 
$or_{X\times X} \otimes or^{-1}_{\Delta_X}$ on the 
standard brane $L_{\Delta_X}$.
\end{proof}


\subsection{Triangulation of diagonal}
We explain here how the choice of a triangulation of $X$ allows us to express the 
diagonal brane $L_{\Delta_X}$ in terms of costandard branes. The primary content of the
section is in developing notation and collecting preliminaries 
for the arguments of subsequent sections.

\medskip

Fix a triangulation $\CT=\{\tau_\alpha\}$ of $X$, and
consider the $\CT$-constructible dg derived category $Sh_\CT(X)$.
Consider the inclusions $j_\fra:\tau_\fra\hookrightarrow X$,
and the corresponding standard sheaves $j_{\fra*} \C_{\tau_\alpha}$.
By Lemma~\ref{standard envelope}, we can express any object of $\Sh_{\CT}(X)$, and
in particular the constant sheaf $\C_X$, as an iterated cone
of maps among the standard sheaves $j_{\fra*} \C_{\tau_\fra}$.

Identify $X$ with the diagonal $\Delta_X\subset X\times X$ (via either projection), and
consider the induced triangulation $\Delta_\CT=\{\Delta_{\tau_\fra}\}$ of $\Delta_X$.
Consider the inclusions $d_\alpha:\Delta_{\tau_\fra}\hookrightarrow \Delta_X$,
and the corresponding standard sheaves $d_{\fra*} \C_{\Delta_{\tau_\fra}}$.
Again, by Lemma~\ref{standard envelope}, we can express
the constant sheaf $\C_{\Delta_X}$ as an iterated cone
of maps among the standard sheaves $d_{\fra*} \C_{\Delta_{\tau_\fra}}$.

\medskip

Next, recall that the diagonal brane $L_{\Delta_X}$ is the standard object of $F(T^*X \times T^*X)$
obtained as the microlocalization
$$
L_{\Delta_X} \simeq \mu_{X\times X}(\C_{\Delta_X}).
$$
By construction, its underlying Lagrangian is the conormal bundle
$
T^*_{\Delta_X} (X\times X).
$

For each simplex of $\CT$ consider the standard object $L_{\Delta_{\tau_\fra}*}$
of $F(T^*X \times T^*X)$
obtained as the microlocalization
$$
L_{\Delta_{\tau_\fra}*} \simeq \mu_{X\times X}(d_{\alpha!}\C_{\Delta_{\tau_\fra}}).
$$
By construction, we can take its underlying Lagrangian to be in the following form.
Fix a non-negative function $m_\fra:X\to \R$ that vanishes precisely on the boundary 
$\del \tau_\fra\subset X$,
and consider the function $f_\fra:{\tau_\fra} \to \R$ given by the logarithm 
$f_\fra = \log m_\fra$.
Then the underlying Lagrangian of $L_{\Delta_{\tau_\fra}!}$ can be identified with
the fiberwise sum
$$
T^*_{\Delta_{\tau_\fra}}(X\times X)   +  \Gamma_{p_2^*d f_\fra} \subset T^*(X\times X),
$$
where $\Gamma_{p_2^*d f_\fra}$ denotes the graph of the pullback
${p_2^*d f_\fra}$ via projection to the second factor 
$$
p_2:X\times X\to X. 
$$

Since we can express the constant sheaf $\C_{\Delta_X}$ as a twisted complex
of standard sheaves $d_{\fra*} \C_{\Delta_{\tau_\fra}}$,
we can express the brane $L_{\Delta_X}$ as a twisted complex
of standard branes $L_{\Delta_{\tau_\fra}*}$.
It will be convenient to have a slightly modified version
of the preceding as recorded in the following.

\begin{lem}\label{lem cut}
The dual brane $\alpha_X(L_{\Delta_X})$ can be expressed as a twisted complex of
the standard branes $L_{\Delta_{\tau_\fra}*}$.
\end{lem}

\begin{proof}
Tracing through the definitions, we have the elementary identity
$\alpha_X(L_{\Delta_X}) =\mu_X(\CD_X(\C_{\Delta_X}))$.
Thus we can repeat the preceding discussion replacing the constant sheaf
$\C_{\Delta_X}$ by its Verdier dual $\CD_X(\C_{\Delta_X})$.
\end{proof}


\subsection{Moving the diagonal} We continue with the notations of the preceding section.

\medskip

Fix a point $x_\fra\in \tau_\fra$, and consider the standard sheaf
$\C_{\{x_\fra\}} \times j_{\fra*} \C_{\tau_\fra}$ as an object of $Sh_c(X\times X)$.

\medskip

Consider the standard object $L_{\{x_\fra\} \times  \tau_\fra *}$ of $F(T^*X \times T^*X)$
obtained as the microlocalization
$$
L_{\{x_\fra\} \times  \tau_\fra *} \simeq 
\mu_{X\times X}(\C_{\{x_\fra\}} \times j_{\fra*} \C_{\tau_\fra}).
$$
By construction, we can take its underlying Lagrangian to be 
the fiberwise sum
$$
T^*_{\{x_\fra\} \times  \tau_\fra}(X\times X)   +  \Gamma_{p_2^*d f_\fra}.
$$
Observe that $L_{\{x_\fra\} \times  \tau_\fra *}$ is the external product
$$
L_{\{x_\fra\} \times  \tau_\fra *} \simeq L_{\{x_\fra\}!} \times L_{\tau_\fra *}
$$
of the factors
$$
L_{\{x_\fra\}} \simeq \mu_X(\C_{\{x_\fra\}})
\qquad
L_{\tau_\fra *} \simeq \mu_X(j_{\fra*} \C_{\tau_\fra}).
$$

\medskip

Consider the conical Lagrangian $\Lambda_\CT\subset T^*X$ given by the 
union of the conormal
bundles of the simplices of the triangulation
$$
\Lambda_\CT = \bigsqcup_\fra T_{\tau_\fra}^* X.
$$
We will use the results of Section~\ref{non-char} to verify the following.

\begin{prop}\label{key point}
For any test objects $P_0, P_1$ of $F(T^*X)$, with 
$P_0^\infty\subset\Lambda_\CT^\infty$,
there is a functorial quasi-isomorphism of complexes
$$
\hom_{F(T^*X\times T^*X)}(P_0\times P_1, L_{\Delta_{\tau_{\fra}}*} )
 \simeq 
\hom_{F(T^*X)}(P_0, L_{\{x_\fra\}}) \otimes \hom_{F(T^*X)}( P_1, L_{\tau_\fra*}).
$$
\end{prop}

\begin{proof}
Since $\tau_\fra$ is contractible, we can find a smooth deformation retract 
$$
\psi_t:\tau_\fra\to \tau_\fra
\qquad
\psi_0 = \id_{\tau_\fra}
\qquad \psi_1(\tau_\fra) = x_\fra.
$$ 
Using $\psi_t$, we define the family of submanifolds 
$$
\tau_{\fra,t} =\{(x_0, x_1) \in X\times X | x_1\in \tau_\fra, x_0 = \psi_t(x_1)\}
$$
satisfying the obvious identifications
$$
\tau_{\fra, 0} = \Delta_{\tau_\fra}
\qquad
\tau_{\fra, 1} = \{x_\fra\} \times \tau_\fra.
$$

Consider the inclusion $d_{\fra, t}: \tau_{\fra, t}\hookrightarrow X\times X$,
the corresponding standard sheaf
$d_{\fra, t*}(\C_{\tau_\fra, t})$, and its microlocalization
$$
L_{\tau_{\fra, t}*} = \mu_{X\times X}(d_{\fra, t*}(\C_{\tau_\fra, t})).
$$
By construction, we can take its underlying Lagrangian to be 
the fiberwise sum
$$
T^*_{\tau_{\fra,t}}(X\times X)   *  \Gamma_{p_2^*d f_\fra}.
$$
We have the obvious identifications
$$
L_{\tau_\fra, 0*} = L_{\Delta_{\tau_\fra}*}
\qquad
L_{\tau_\fra, 1*} = L_{\{x_\fra\} \times  \tau_\fra *}
= L_{\{x_\fra\}} \times L_{\tau_\fra *}.
$$

Now fix test objects $P_0,P_1$ of $F(T^*X)$,
with $P_0^\infty\subset\Lambda_\CT^\infty$, and let $\Lambda\subset T^*X$
be a conical Lagrangian with $P_1^\infty \subset \Lambda^\infty$.
We would like to show the family of Floer complexes
$$
\hom_{F(T^*X\times T^*X)}(P_0\times P_1, L_{\Delta_{\tau_{\fra,t}}*})
$$
has constant cohomology with respect to $t$.
Choose a small $\eta >0$,
and consider the submanifold $\tau_{\fra,\eta} = \{x\in \tau_\fra | m_\fra(x)>\eta\}$,
and
the family of submanifolds 
$$
\tau_{\fra,\eta, t} =\{(x_0, x_1) \in X\times X | x_1\in \tau_{\fra,\eta}, x_0 = \psi_t(x_1)\}.
$$
Consider the inclusion $d_{\fra, \eta, t}: \tau_{\fra, \eta, t}\hookrightarrow X\times X$,
the corresponding costandard sheaf
$d_{\fra, \eta, t!}(\C_{\tau_\fra, \eta, t})$, and its microlocalization
$$
L_{\tau_\fra, \eta, t*} = \mu_{X\times X}(d_{\fra, \eta, t*}(\C_{\tau_\fra, \eta, t})).
$$
By construction, we can take its underlying Lagrangian to be the fiberwise sum
$$
T^*_{\Delta_{\tau_\fra}}(X\times X)   +  \Gamma_{p_2^*d f_{\fra,\eta}}
$$
where $f_{\fra,\eta}:\tau_{\fra,\eta}\to \R$ is the function given by the logarithm 
$f_{\fra,\eta} = \log(m_\fra -\eta)$.
Then for sufficiently small $\eta>0$, and all $t$, 
we have a quasi-isomorphism
$$
\hom_{F(T^*X\times T^*X)}(P_0\times P_1, L_{\Delta_{\tau_{\fra, t}}*} )
\simeq
\hom_{F(T^*X\times T^*X)}(P_0\times P_1, L_{\Delta_{T_{\fra, \eta, t}}*}).
$$

Finally, by construction, the family of branes $L_{\tau_\fra, \eta, t*}$ is non-characteristic
with respect to the product conical Lagrangian $\Lambda_\CT\times \Lambda$.
Thus by Proposition~\ref{prop brane non-char} (see in particular Example~\ref{ex brane non-char}),
the family of complexes 
$$
\hom_{F(T^*X\times T^*X)}(P_0\times P_1, L_{\Delta_{T_{\fra, \eta, t}}*})
$$
has constant cohomology.
Furthermore, at $t=0$, the family calculates the left hand side of the proposition, and at $t=1$,
it calculates the right hand side (cf. Lemma~\ref{lem product floer calc}).
\end{proof}


\subsection{Proof of Theorem~\ref{theorem projection}} 
Now let us wrap up the proof of Theorem~\ref{theorem projection}.
We continue with the notation of the preceding sections.

\medskip

Fix once and for all an object $L$ of $F(T^*X)$.
Our aim is to show that for any test object $P$ of $F(T^*X)$,
we can functorially express the Floer complex 
$\hom_{F(T^*X)}(P,L)$ as a twisted complex of the Floer complexes
$\hom_{F(T^*X)}(P,L_{Y*})$ of standard branes $L_{Y*}\hra T^*X$ associated
to submanifolds $Y\hookrightarrow X$.

\medskip

First, by Proposition~\ref{diag}, 
there is a functorial quasi-isomorphism of Floer complexes
$$
\hom_{F(T^*X)}(P, L)\simeq \hom_{F(T^*X\times T^*X)}(L_{\Delta_X}, L \times\alpha_X(P)).
$$
It will be convenient to rewrite the preceding
in a slightly modified form.
Since the brane duality $\alpha_X$ is an anti-equivalence, we have a functorial
quasi-isomorphism
$$
\hom_{F(T^*X)}(P, L)\simeq \hom_{F(T^*X\times T^*X)}(\alpha_X(L) \times P,\alpha_X(L_{\Delta_X})).
$$

Next, fix a triangulation $\CT=\{\tau_\alpha\}$ of $X$ along with the induced
triangulation $\Delta_{\CT}=\{\Delta_{\tau_\alpha}\}$ of $\Delta_X\subset X\times X$.
By Lemma~\ref{lem cut},
we can express $\alpha_X(L_{\Delta_X})$ as a twisted complex of the standard branes
$L_{\Delta_{\tau_\fra}*}$.

Suppose further that $\CT$ is chosen fine enough so that
$
L^\infty \subset \Lambda_\CT^\infty.
$
Then by Proposition~\ref{key point},
there is a functorial quasi-isomorphism of complexes
$$
\hom_{F(T^*X\times T^*X)}(\alpha_X(L)\times P, L_{\Delta_{\tau_{\fra}}*})
 \simeq 
\hom_{F(T^*X)}(\alpha_X(L), L_{\{x_\fra\}}) \otimes \hom_{F(T^*X)}(P, L_{\tau_\fra*}).
$$

Putting together the preceding identifications,
we have functorially expressed the Floer complex $\hom_{F(T^*X)}(P,L)$ as a 
twisted complex with terms the Floer complexes 
$\hom_{F(T^*X)}( \alpha_X(P), L_{\tau_\fra*})$.
This completes the proof of Theorem~\ref{theorem projection}.
\hfill$\square$

\begin{rmk}\label{rem coeffs}
In the expression of the brane $L$ as a twisted complex of the standard branes
$L_{\tau_\fra*},$
the coefficients appearing are the functionals 
$$
\hom_{F(T^*X)}(\alpha_X(L), L_{\{x_\fra\}}) \simeq \hom_{F(T^*X)}(L_{\{x_\fra\}}, L).
$$

Although we will not use it, the proof of Theorem~\ref{theorem projection}
leads to the following precise form of $L$ as a twisted complex.
First, we can identify the coefficients $\hom_{F(T^*X)}(L_{\{x_\fra\}}, L)$ with the 
shifted Floer complexes $\hom_{F(T^*X)}(L_{\tau_\fra!}, L)$.
Next, the triangulation provides the structure of dual cell complexes on the collection of
standard branes $L_{\tau_\alpha*}$ and costandard
branes $L_{\tau_\fra!}$ associated to the simplices $\tau_\fra$.
Finally, this induces the structure of a twisted complex on the branes
$ \hom_{F(T^*X)}(L_{\tau_\fra!}, L)\otimes L_{\tau_\fra*}$ appearing
in the decomposition of $L$. 
\end{rmk}


\subsection{From branes to sheaves}\label{sect from branes to sheaves}

Since microlocalization is a quasi-equivalence
$$
\xymatrix{
\mu_X: Sh_c(X) \ar[r]^-{\sim} & F(T^*X),
}$$
the corresponding pullback of right $A_\infty$-modules is also
a quasi-equivalence 
$$
\xymatrix{
\mu_X^*: \rmod(F(T^*X)) \ar[r]^-{\sim} & \rmod(Sh_c(X)).
}$$
In particular, given an object $L$ of $F(T^*X)$, we can take its Yoneda
module $\CY_r(L)$, and ask what object $\CF$ of $Sh_c(X)$ 
quasi-represents $\mu_X^*\CY_r(L)$.

Here is an informal way to think about an object $\CF$ 
that quasi-represents $\mu_X^*\CY_r(L)$.
Given an open submanifold $i:U\hra X$, we have quasi-isomorphisms
of complexes
$$
\CF(U) \simeq \hom_{Sh_c(X)}(i_!\C_U, \CF)
\simeq \hom_{F(T^*X)}(L_{U!}\otimes or_X[-\dim X], L).
$$
The first quasi-isomorphism is by adjunction, and the second is by the fact that 
$\CF$ 
quasi-represents $\mu_X^*\CY_r(L)$.
Here we have taken $L_{U!}\otimes or_X[-\dim X]$
since $\omega_X \simeq or_X[\dim X]$. 

Similarly, given an inclusion of open submanifolds $i_{0}^1:U_0\hra U_1$, 
we have a diagram that commutes at the level of cohomology
$$
\xymatrix{
\CF(U_1) \ar[d]_-{\wr}\ar[r]^-{\rho^1_0} & \CF(U_0) \ar[d]^-{\wr}\\
\hom_{F(T^*X)}(L_{U_1!}\otimes or^\vee_X[-\dim X], L)
\ar[r]^-{\iota_0^1 } & 
 \hom_{F(T^*X)}(L_{U_0!}\otimes or^\vee_X[-\dim X], L).
}$$
Here $\rho^1_0$ denotes the sheaf restriction map,
and $\iota_0^1$ denotes the map induced by the canonical degree zero morphism in
the complex
$$
\hom_{F(T^*X)}(L_{U_0!}, L_{U_1!}) \simeq \hom_{Sh_c(X)}(i_{0!}\C_{U_0}, i_{1!}\C_{U_1})
\simeq (\Omega(U_0), d).
$$

Finally, we record the following consequence of the proof of Theorem~\ref{theorem projection}.
Fix a conical Lagrangian $\Lambda\hra T^*X$,
and let $F(T^*X)_\Lambda$ be the full subcategory of $F(T^*X)$ of twisted complexes
of Lagrangian branes $L$ such that $L^\infty\subset \Lambda^\infty$.

\begin{prop}\label{prop constructibility of branes}
Given an object $L$ of $F(T^*X)_\Lambda$, consider its Yoneda
module $\CY_r(L)$, and an object $\CF$ of $Sh_c(X)$ that
quasi-represents the pullback $\mu_X^*\CY_r(L)$.

Then for any stratification $\CS=\{S_\alpha\}$ 
of $X$ such that $\Lambda \subset \Lambda_\CS = \cup_\alpha T^*_{S_\alpha} X$,
the object $\CF$ lies in $Sh_{\CS}(X)$.
\end{prop}

\begin{proof}
In the proof of Theorem~\ref{theorem projection},
we showed that $L$ can be expressed as a twisted complex
of the standard branes $L_{\tau_\frb*}$,
for any triangulation $\CT=\{\tau_\frb\}$ refining $\CS$.
In other words, $\CF$ can be expressed as a twisted complex of the standard sheaves
$j_{\frb*}\C_{\tau_\frb}$.

In fact, in place of the triangulation $\CT$,
we can take
any disjoint cell decomposition $\CC=\{c_\frb\}$ of $X$
refining $\CS$ in the sense
that each cell $j_\frb:c_\frb\hra X$ lies in some stratum $S_\alpha$. 
To see this level of generality, observe that all we need
for the proof of Theorem~\ref{theorem projection}
 is that the dualizing
complex $\omega_X=\CD_X(\C_X)$ can be expressed as a twisted complex of
the standard sheaves $j_{\frb*}\C_{c_\frb}$,
and that each cell $c_\frb$ can be deformation 
contracted within the stratum $S_\alpha$  containing $c_\frb$
to a point $x_\frb\in c_\frb$.

Thus $\CF$ can be expressed as a twisted complex of the standard
sheaves $j_{\frb*}\C_{c_\frb}$ on the cells of any such cell decomposition $\CC=\{c_\frb\}$.
The proposition now follows from Lemma~\ref{constructibility lemma}
immediately below.
%
\end{proof}

%


\begin{lem}\label{constructibility lemma}
Let $\CF$ be an object of $Sh_c(X)$. Suppose that for any 
 cell decomposition $\CC=\{c_\frb\}$ of $X$ refining the stratification $\CS=\{S_\alpha\}$,
we can express $\CF$ as a twisted complex of the standard sheaves $j_{\frb*}\C_{c_\frb}$
on the cells.
Then  $\CF$ is $\CS$-constructible.
\end{lem}

\begin{proof}
Fix any point $p \in X$, and let $S_\alpha$ be the stratum of $\CS$ containing $p$. 
Choose an open ball $B_{\alpha, p}\subset S_\alpha$ containing $p$, and a normal slice $N_{\alpha, p} \subset X$ to $S_\alpha$ at $p$.
We will consider $B_{\alpha, p}$ as a smooth manifold, and equip $N_{\alpha, p}$ with the stratification induced by restricting $\CS$.

By the Thom Isotopy Lemma, there is a stratum preserving
homeomorphism from the product $B_{\alpha, p}\times N_{\alpha, p}$,
equipped with the product stratification, to a neighborhood $U_p
\subset X$
of the point $p$, with the stratification induced by restricting $\CS$. Furthermore, the restriction of the homeomorphism to each stratum
is a diffeomorphism.

Choose any triangulation $\CT_{N_{\alpha, p}}$ of the normal slice $N_{\alpha, p}$ refining the stratification 
induced by 
$\CS$.
Consider the induced product stratification
of the neighborhood
$U_p \simeq B_{\alpha, p}\times N_{\alpha, p}$. 
Extend this to any cell decomposition $\CC=\{c_\frb\}$ of all of $X$
by cutting up the complement $X\setminus U_p$ into cells.

By construction, the restriction of the standard sheaf $j_{\frb*}\C_{c_\frb}$
of any cell $c_\frb$ of $\CC$
to the cell $B_{\alpha, p}$ is constant. Hence by assumption,
the restriction of $\CF$ to the cell $B_{\alpha, p}$
is constant as well. Thus the restriction of $\CF$ to the stratum $S_\alpha$ is locally constant. 
\end{proof}


\section{Functoriality}

Now that we have established that microlocalization is a quasi-equivalence
$$
\xymatrix{
\mu_X: Sh_c(X) \ar[r]^-{\sim} & F(T^*X),
}
$$
we can collect some formal consequences for future applications.


\subsection{Duality revisited}\label{sec duality revisited}

Recall the brane duality equivalence
$$
\alpha_{X}:F(T^*X)^\circ \risom F(T^*X)
$$
introduced in Section~\ref{sec duality}.
Our aim in this section is to confirm that microlocalization $\mu_X$ intertwines brane duality
with Verdier
duality 
$$
\CD_{X}:Sh_c(X)^\circ \risom Sh_c(X).
$$

For a submanifold $i_Y:Y\hra X$, Verdier duality exchanges the associated standard
and costandard objects
$$
i_{Y*}\C_Y \simeq \CD(i_{Y!}\omega_Y).
$$
Likewise, by construction, brane duality exchanges 
the associated standard
and costandard branes
$$
\alpha_X(L_{Y*})\simeq L_{Y!}.
$$
%

\begin{prop}\label{dual confirm}
There is a quasi-isomorphism
$$
\mu_X\circ\CD_X \simeq \alpha_X\circ \mu_X: \Sh_c(X)^\circ \to F(T^*X).
$$
\end{prop}

\begin{proof}
By Proposition~\ref{embed costandards},
for any open submanifold $i_U:U\hra X$, there is a quasi-isomorphism of left $\Sh_c(X)$-modules
$$
\mu_X^*(\CY_{\ell}(L_{U!})) \simeq \CY_\ell(i_{U!}\omega_U):Sh_c(X)\to \Ch.
$$
In other words, for any test object $\CF$ of $Sh_c(X)$, there is a functorial quasi-isomorphism
$$
\hom_{F(T^*X)}(L_{U!}, \mu_X(\CF))\simeq
\hom_{Sh_c(X)}(i_{U!}\omega_U, \CF)
$$
Since microlocalization $\mu_X$ is a quasi-embedding, there is a functorial quasi-isomorphism
$$
\hom_{Sh_c(X)}(i_{U!}\omega_U, \CF)
\simeq
\hom_{F(T^*X)}(\mu_X(i_{U!}\omega_U), \mu_X(\CF)).
$$
Finally, since $\mu_X$ is a quasi-equivalence (so objects of the form $\mu_X(\CF)$
generate $F(T^*X)$), and the Yoneda functor
$\CY_\ell$ is a quasi-embedding, there is an isomorphism
$$
\mu_X(i_{U!}\omega_U) \simeq L_{U!}.
$$

We conclude that there are isomorphisms
$$
\mu_X(\CD_X(i_{U*}\C_U)) \simeq
\mu_X(i_{U!}\omega_U)
\simeq
L_{U!}
 \simeq \alpha_X(L_{U*})
 \simeq
 \alpha_X(\mu_X(i_{U*}\C_U)).
$$
By Lemma~\ref{standard star envelope}, standard objects $i_{U*}\C_U$ generate $Sh_c(X)$,
and so we have the asserted quasi-isomorphism of functors.
\end{proof}


\subsection{Integral transforms for sheaves }\label{transform}
In this and the following section, we describe standard functors for sheaves and branes.
At this stage, their compatibility is completely formal: 
it depends only on the fact that microlocalization
is a quasi-equivalence intertwining Verdier duality with brane duality.

\medskip

Given two real analytic manifolds $X_0, X_1$, consider the standard projections 
$$
p_{0}:X_0\times X_1\to X_0
\qquad 
p_{1}:X_0\times X_1\to X_1.
$$

\medskip

Following~\cite[Section 3.6]{KS}, for an object $\CK$ of $Sh_c(X_0\times X_1)$, we have the integral transforms
$$
\Phi_\CK^*: Sh_c(X_1) \to Sh_c(X_0)
\qquad
\Phi_{\CK*}: Sh_c(X_0)  \to Sh_c(X_1)
$$
$$
\Phi_\CK^*(\CF) = p_{0!}(\CK\otimes p_{1}^*(\CF))
\qquad
\Phi_{\CK*}(\CF) = p_{1*}(\intHom(\CK, p_0^!(\CF)))
$$
Similarly, reversing the roles of $X_0$ and $X_1$, we have the integral transforms
$$
\Phi_{\CK!}: Sh_c(X_0)  \to Sh_c(X_1)
\qquad 
\Phi_\CK^!: Sh_c(X_1) \to Sh_c(X_0)
$$
$$
\Phi_{\CK!}(\CF) = p_{1!}(\CK \otimes p_0^*(\CF)))
\qquad
\Phi_\CK^!(\CF) = p_{0*}(\intHom(\CK, p_1^!(\CF)))
$$

Standard identities 
imply the following:
(1) $(\Phi^*_\CK,\Phi_{\CK*})$ and $(\Phi_{\CK!},\Phi^!_{\CK})$ are each adjoint pairs.
(2) The construction is functorial in $\CK$ in the sense that we have functionals
$$
\Phi^*\simeq\Phi_{*}: Sh_c(X_0) \otimes Sh_c(X_0\times X_1)^\circ \otimes Sh_c(X_1)^\circ \to \Ch
$$
$$
(\CF_1, \CK, \CF_0) \mapsto \hom_{Sh_c(X_0)}(\Phi_\CK^*(\CF_1), \CF_0) \simeq \hom_{Sh_c(X_1)}(\CF_1, \Phi_{\CK*}(\CF_0))
$$
$$
\Phi_!\simeq\Phi^{!}: Sh_c(X_0)^\circ \otimes Sh_c(X_0\times X_1)^\circ \otimes Sh_c(X_1) \to \Ch
$$
$$
(\CF_0, \CK, \CF_1) \mapsto \hom_{Sh_c(X_1)}(\Phi_{\CK!}(\CF_0), \CF_1) \simeq \hom_{Sh_c(X_0)}(\CF_0, \Phi_{\CK}^!(\CF_1))
$$
(3) We have quasi-isomorphisms of functors
$$
\Phi^*_{\CK} \simeq \CD_{X_0}\circ \Phi^!_{\CK} \circ\CD_{X_1}
\qquad
\Phi_{\CK*} \simeq \CD_{X_1}\circ \Phi_{\CK!} \circ\CD_{X_0}
$$
$$
\Phi_{\CK!} \simeq \CD_{X_1}\circ \Phi_{\CK*} \circ\CD_{X_0}
\qquad
\Phi^!_{\CK} \simeq \CD_{X_0}\circ \Phi^*_{\CK} \circ\CD_{X_1}
$$

\begin{ex}
Our notation is motivated by the following example. 
Fix a $\CC$-map $\ff:X_0\to X_1$,
consider the graph 
$$
\Gamma_\ff=\{ (x_0, x_1) \in X_0\times X_1| x_1 =\ff(x_0)\},
$$
and let $\C_{\Gamma_\ff}$ denote the constant sheaf along $\Gamma_\ff$.
Then we have canonical identifications of functors
$$
(\ff^*, \ff_*) \simeq (\Phi_{\C_{\Gamma_\ff}}^*,\Phi_{\C_{\Gamma_\ff}*})
\qquad
(\ff_!, \ff^!) \simeq (\Phi_{\C_{\Gamma_\ff}!},\Phi_{\C_{\Gamma_\ff}}^!)
$$
\end{ex}


\subsection{Integral transforms for branes}
We discuss here the analogous integral transforms associated to 
objects of $F(T^*X_0\times T^*X_1)$. 

For any object $L$ of $F(T^*X_0\times T^*X_1)$, we have functors
$$
\tilde \Psi_L^*: F(T^*X_1) \to  \lmod(F(T^*X_0))^\circ
$$
$$
\tilde \Psi_L^*(P_1): P_0\mapsto \hom_{F(T^*X_0\times T^*X_1)}(L, P_0\times \alpha_{X_1}(P_1))
$$
$$
\tilde\Psi_{L*}: F(T^*X_0) \to  \rmod(F(T^*X_1))
$$
$$
\tilde\Psi_{L*}(P_0) :P_1\mapsto  \hom_{F(T^*X_0\times T^*X_1)}(L, P_0 \times \alpha_{X_1}(P_1))
$$
$$
\tilde\Psi_{L!}: F(T^*X_0) \to  \lmod(F(T^*X_1))^\circ
$$
$$
\tilde\Psi_{L!}(P_0) :P_1\mapsto \hom_{F(T^*X_0\times T^*X_1)}(L, \alpha_{X_0}(P_0)\times P_1)
$$
$$
\tilde\Psi_{L}^!: F(T^*X_1) \to  \rmod(F(T^*X_0))
$$
$$
\tilde\Psi_{L}^!(P_1) :P_0 \mapsto  \hom_{F(T^*X_0\times T^*X_1)}(L, \alpha_{X_0}(P_0)\times P_1)
$$
Note that the constructions are
functorial in $L$ in the contravariant sense.

\begin{prop}
Consider an object $\CK$ of $Sh_c(X_0\times X_1)$,
and its microlocalization $L = \mu_{X_0\times X_1}(\CK)$.
Then there are functorial quasi-isomorphisms 
$$
\CY_\ell\circ \mu_{X_0}\circ\Phi_\CK^*
\simeq
\tilde \Psi_{L}^*\circ\mu_{X_1}
\qquad
\CY_r\circ\mu_{X_1}\circ\Phi_{\CK*}
\simeq
\tilde\Psi_{L*}\circ\mu_{X_0}
$$
$$
\CY_\ell\circ\mu_{X_1}\circ\Phi_{\CK !}
\simeq
\tilde\Psi_{L!}\circ\mu_{X_0}
\qquad
\CY_r\circ\mu_{X_0}\circ\Phi^!_{\CK}
\simeq
\tilde\Psi^!_{L}\circ\mu_{X_1}
$$
\end{prop}

\begin{proof}
We establish the second quasi-isomorphism (the case of the usual pushforward); the arguments for the others are similar.
It suffices to consider test objects $L_0$ of $F(T^*X_0)$ and $L_1$ of $F(T^*X_1)$
of the form $L_0=\mu_{X_1}(\CF_0)$ and $L_1= \mu_{X_1}(\CF_1)$,
and to establish a functorial quasi-isomorphism
$$
\hom_{\Sh_c(X_0\times X_1)}(\CF_1, p_{1*}(\intHom(\CK, p_0^!(\CF))))
\simeq
\hom_{F(T^*X_0\times T^*X_1)}(L, L_0 \times \alpha_{X_1}(L_1)).
$$
By standard identities, this is nothing more than a functorial quasi-isomorphism
$$
\hom_{\Sh_c(X_0\times X_1)}(\CK, p^*_0(\CF_0) \otimes p_1^*(\CD_{X_1}(\CF_1)))
\simeq
\hom_{F(T^*X_0\times T^*X_1)}(L, L_0 \times \alpha_{X_1}(L_1)).
$$
Now the assertion follows immediately from Proposition~\ref{dual confirm}.
\end{proof}

We see from the proposition that the modules arising from the functors
$
\tilde \Psi_{L}^*,
$
$
\tilde\Psi_{L*},
$
$\tilde\Psi_{L!},
$
$\tilde\Psi^!_{L}$
are representable.
Namely, for $L\simeq \mu_{X_0\times X_1}(\CK)$,
we can take the representing functors to be the compositions
$$
\Psi_L^* = \mu_{X_0} \circ\Phi^*_\CK \circ \pi_{X_1}
\qquad
\Psi_{L*} = \mu_{X_1} \circ\Phi_{\CK*} \circ \pi_{X_0}
$$
$$
\Psi_{L!} = \mu_{X_1} \circ\Phi_{\CK!} \circ \pi_{X_0}
\qquad
\Psi_{L}^! = \mu_{X_0} \circ\Phi_{\CK}^! \circ \pi_{X_1}
$$

Their basic properties are immediate from the definitions:
(1) $(\Psi_{L}^*, \Psi_{L*})$ and $(\Psi_{L!}, \Psi_L^!)$ are each adjoint pairs.
(2) The construction is functorial in $L$ in the sense that we have functors
$$
\Phi^*\simeq\Phi_{*}: F(T^*X_0\times T^*X_1) 
\to \mbox{{\em fun}}_{A_\infty}(F(T^*X_0), \rmod(F(T^*X_1))
$$
$$
\Psi_!\simeq\Psi^{!}: F(T^*X_0\times T^*X_1) 
\to \mbox{{\em fun}}_{A_\infty}(F(T^*X_1), \rmod(F(T^*X_0))
$$
where $\mbox{{\em fun}}_{A_\infty}$ denotes the $A_\infty$-category of $A_\infty$-functors.
(We have written the functors in this form rather than as functionals
since the notion of internal hom for $A_\infty$-categories
is straightforward, while that of tensor product is more delicate.)
(3) There are functorial quasi-isomorphisms
$$\
\Psi_{L}^*\simeq \alpha_{X_0}\circ\Psi_{L}^!\circ\alpha_{X_1}
\qquad
\Psi_{L*}\simeq \alpha_{X_1}\circ\Psi_{L!}\circ\alpha_{X_0}
$$
$$
\Psi_{L!}\simeq \alpha_{X_1}\circ\Psi_{L*}\circ\alpha_{X_0}
\qquad
\Psi_{L}^!\simeq \alpha_{X_0}\circ\Psi_{L}^*\circ\alpha_{X_1}
$$
intertwining the functors with brane duality.

\begin{ex}

Let $\ff:X_0\to X_1$ be a $\CC$-map.
Consider the graph 
$$
\Gamma_\ff=\{ (x_0, x_1) \in X_0\times X_1| x_1 =\ff(x_0)\},
$$
and let $\C_{\Gamma_\ff}$ denote the constant sheaf along $\Gamma_\ff$.

Let $L_{\ff}$ be the standard object of $F(T^*X_0 \times T^*X_1)$
obtained as the microlocalization
$$
L_{\ff} \simeq \mu_{X_0\times X_1}(\C_{\Gamma_\ff}).
$$
By construction, when $\ff$ is smooth, we can take the Lagrangian underlying $L_{\ff}$ 
to be the conormal bundle
$
T^*_{\Gamma_\ff} (X_0\times X_1).
$

Applying the above constructions, we obtain functors 
$\Psi_{L_\ff }^*,\Psi_{L_\ff *}, \Psi_{L_\ff !},\Psi_{L_\ff }^!.$

\end{ex}

\begin{cor}
For any $\CC$-map $\ff:X_0\to X_1$,
there are quasi-isomorphisms 
$$
\Psi_{L_\ff}^*\circ \mu_{X_1} \simeq \mu_{X_0} \circ \ff^*
\qquad
\Psi_{L_\ff*}\circ \mu_{X_0} \simeq \mu_{X_1} \circ \ff_*
$$
$$
\Psi_{L_\ff!}\circ \mu_{X_0} \simeq \mu_{X_1} \circ \ff_!
\qquad
\Psi_{L_\ff}^!\circ \mu_{X_1} \simeq \mu_{X_0} \circ \ff^!
$$
\end{cor}


\subsection{Correspondence interpretation}\label{sec corr}

In this informal section, we sketch how the integral transforms of the preceding sections
are related to the beautiful theory of quilted Riemann surfaces and generalized branes 
mathematically developed by Wehrheim-Woodward \cite{WW} (and 
pioneered from a physical perspective by Khovanov-Rozansky~\cite{KR}
under the name world-sheet foam).
We do not use the discussion of this section and include it for the interested reader already 
familiar with the constructions of~\cite{WW}.

In what follows, we assume that all of our manifolds are orientable, so that their
cotangent bundles are spin. The main reason for imposing this condition will be
that for a product $T^*X_0\times T^*X_1$,
the canonical background class 
will then be the product of the canonical background classes
$$
(\pi_0\times \pi_1)^*w_2(X_0 \times X_1)
=
\pi_0^*w_2(X_0) + \pi_1^*w_2(X_1).
$$

\subsubsection{Generalized branes}
By a generalized Lagrangian submanifold of $T^*X$, we mean
a sequence of compact real analytic manifolds $pt = X_{-m},$ $X_{-m+1},
\ldots, X_{-1}, X_0 = X$,
for some $m>0$, and a sequence 
of Lagrangian submanifolds $\ul L= (L_{(-m,-m+1)}, \ldots, L_{(-1,0)})$ 
in the successive products
$$
L_{(-k,-k+1)}\subset (T^*X_{-k})^-\times T^*X_{-k+1},
\qquad
\mbox{ for $k=1,\ldots, m$}.
$$
As usual, to control the behavior of $L_{(-k,-k+1)}$ near infinity, 
we require that
its closure
in the product compactification 
is a $\CC$-subset, and that there is a perturbation
to a tame Lagrangian.

A brane structure on a generalized Lagrangian submanifold $\ul L$
consists of a sequence $\ul \CE= (\CE_{(-m,-m+1)}, \ldots, \CE_{(-1,0)})$ of flat vector
bundles
$$
\CE_{(-k,-k+1)}\to L_{(-k,-k+1)},
\qquad
\mbox{ for $k=1,\ldots, m$}.
$$
and a sequence 
of gradings
and relative pin structures. 
For simplicity,
we will take the gradings and relative pin structures to be defined with respect to the canonical
product bicanonical trivializations
and background forms
respectively.

\subsubsection{Composition of correspondences}

Following Wehrheim-Woodward~\cite{WW},
there is
a triangulated category $DF^\#(T^*X)$ 
whose objects are twisted complexes of 
generalized Lagrangian branes. 
Work in progress of Mau-Wehrheim-Woodward~\cite{MWW} will provide
an $A_\infty$-enhancement of this story but we content ourselves
here with discussing things at the cohomological level.

A primary motivation for introducing generalized Lagrangian branes is that 
Lagrangian correspondences act on them: there is
a triangulated functor
$$
DF^\#(T^*X_0) \otimes DF((T^*X_0)^-\times T^*X_1) \to DF^\#(T^*X_1)
$$
given on objects by concatenation
$$
(\ul L=(L_{(-m,-m+1)}, \ldots, L_{(-1,0)}), L_{(0,1)}) \mapsto 
\ul L\# L_{(0,1)}= (L_{(-m,-m+1)}, \ldots, L_{(-1,0)}, L_{(0,1)}).
$$
The structure of the categories $DF^\#(T^*X_0), DF^\#(T^*X_1)$ and the 
composition functor are given by counting quilted Riemann surfaces.
In particular, there is a Floer functional to chain complexes
$$
DF^\#(T^*X_0)^\circ \otimes DF((T^*X_0)^-\otimes T^*X_1)^\circ \otimes DF^\#(T^*X_1) \to D(\Ch)
$$
$$
(\ul L_0, L_{(0,1)}, \ul L_1)\mapsto 
\hom_{DF^\#(T^*X_1)}(\ul L_0\# L_{(0,1)}, \ul L_1).
$$

Observe that there is an obvious functor $DF(T^*X)\to DF^\#(T^*X)$,
and thus given an object $\ul L$ of $DF^\#(T^*X)$, we can think of it as defining a 
left $DF(T^*X)$-module via the Yoneda map
$$
\CY_{\ell,ord}:DF^\#(T^*X)\to\lmod(DF(T^*X))
$$
$$
\CY_{\ell,ord}(\ul L): P \mapsto \hom_{F^\#(T^*X)}(\ul L, P)
$$


\subsubsection{Compatibility}

Now consider the product symplectomorphism
$$
a_0\times\id_1:T^* X_0 \times T^*X_1 \risom (T^*X_0)^-\times T^*X_1
$$
$$
(a_0\times\id_1)(x_0,\xi_0; x_1, \xi_1)= (x_0,-\xi_0; x_1, \xi_1).
$$
It induces an equivalence by transport of structure
$$
(a_0\times\id_1)_*:DF(T^*X_0 \times T^*X_1) \risom DF((T^*X_0)^-\times T^*X_1).
$$

To ensure the compatibility of the following proposition, we introduce
a twisted version of the above equivalence. Namely, we define the equivalence
$$
(a_0\times\id_1)^\sim_*:DF(T^*X_0 \times T^*X_1) \risom DF((T^*X_0)^-\times T^*X_1)
$$
to be the composition of $(a_0\times\id_1)_*$ with the twist by the pullback $p_0^*(or_{X_0})$
of the orientation bundle from the first factor. 

Given an object $L$ of $DF(T^*X_0 \times T^*X_1)$, we write 
$$
L_{(0,1)}= (a_0\times\id_1)^\sim_*(L)
$$ 
for the corresponding object of $DF((T^*X_0)^-\times T^*X_1)$.

We leave the proof of the following to the interested reader; it is not used in 
other parts of the paper.
Similar identities exist for the other ``integral transforms".

\begin{prop}
Given objects $P_0$ of $DF(T^*X_0)$ and $L$ of $DF(T^*X_0 \times T^*X_1)$,
there is a functorial isomorphism of left $DF(T^*X_1)$-modules
$$
\CY_{\ell, ord}(P_0\# L_{(0,1)}) \simeq \tilde \Psi_{L!}(P_0).
$$
\end{prop}


\section{Appendix: invariance of calculations}\label{appendix}

We discuss here some aspects of the invariance of calculations among microlocal branes.
We assume the standard (though highly intricate) theory for compact exact branes
(in the form explained by Seidel~\cite{Seidel}), and comment about the modest modifications
needed to treat the noncompact branes we consider. We do not attempt
anything approximating a comprehensive discussion, but rather specifically argue for the
independence of the taming perturbation in the definition of a microlocal brane (see
Section~\ref{sect Fukaya category}).


\subsection{Almost complex structures}

Recall that in the definition of $F(T^*X)$,
we work with 
an asymptotically conical almost complex structures $J_{con}\in \on{End}(T (T^*X))$
(see~Section~\ref{sect prels}).
Then we require that every microlocal brane $L$ comes equipped with a taming perturbation
$\psi$ that moves it to a brane $\psi(L)$ that is tame with respect
to the induced metric $g_{con}(v,v)=\omega(v, J_{con}v)$ (see~Section~\ref{sect Fukaya category}). 
These requirements ensure that the moduli spaces defining
the structure constants of the $A_\infty$-operations of $F(T^*X)$ are compact
(see again Section~\ref{sect Fukaya category}). 

Our aim here is to show that in fact calculations among microlocal branes are independent
of the class of asymptotically conical almost complex structures.
For any finite calculation 
(finite number of objects,
finite number of $A_\infty$-operations), we will show that as long as we choose
a compatible almost complex structure $J$ such that $T^*X$ and 
the branes under consideration are tame, the resulting $A_\infty$-operations 
will be compatible with those defined with respect to any other such
almost complex structure $J'$ (in particular, 
an asymptotically conical almost complex structure). 
Furthermore, our arguments can be made compatibly for increasing unions
of finite calculations.

To isolate the role of the almost complex structure, let us fix a finite collection
of branes $L_0, \cdots, L_d \subset T^*X$, and without loss of generality,
assume that they are already 
mutually transverse
and do not intersect each other at infinity.

\begin{lem} \label{indep of ac str}
For any compatible almost complex structures $J, J'$ such that $T^*X$
and the branes $L_0,\ldots, L_d$ are tame, 
there is a quasi-isomorphism between the $A_\infty$-operations
within the collection $L_0,\ldots, L_d$ defined with respect to $J$ and $J'$.
\end{lem}

\begin{proof}
If $J$ and $J'$ coincide outside of a compact set, we refer the reader
to standard homotopy arguments~\cite{Seidel} to construct the sought-after quasi-isomorphism.
So our aim here is to show that we can put ourselves into this situation.

Fix a finite collection of $A_\infty$-operations. Consider the corresponding moduli
problems of $J$-holomorphic disks that calculate the structure constants
of the operations. Recall that we have an a priori diameter bound
on the relevant $J$-holomorphic disks: there is a large $r>0$ such that
none of the disks enter the region of $T^*X$ given by $|\xi|>r/2$.

Now replace $J$ by a compatible almost complex structure $J_{cut}$ of the following form:
\begin{enumerate}
\item $J_{cut}= J$ in the region $|\xi| < r/2$.\
\item $J_{cut}=J'$ in the region $|\xi| > r$.
\end{enumerate}
So $J_{cut}$ equals $J$ in a compact region, $J'$ near infinity,
and whatever one chooses in between.
By construction, for the fixed collection of $A_\infty$-operations, 
we have the same a priori diameter bound
on the relevant $J_{cut}$-holomorphic disks.
This follows from the property (1) above and the local derivation of the diameter bound.
Thus the corresponding moduli spaces for $J$ and for $J_{cut}$ are in fact equal.

Finally, choose a $[0,1]$-family of compatible almost complex structure $J_{t}$
satisfying:
\begin{enumerate}
\item $J_{0}= J_{cut}$.\
\item $J_1 = J'$.\
\item $J_{t}= J_{cut} =J'$ in the region $|\xi| > r$ for all $t$.
\end{enumerate}
Since $J_t$ is constant near infinity, we can apply standard 
homotopy arguments~\cite{Seidel} to compare
the $A_\infty$-operations for $J_{cut}$
and for $J'$. 
\end{proof}


\subsection{Taming pertubations}

We apply here Lemma~\ref{indep of ac str} to see that the choice of taming perturbation $\psi$
in the definition (see~Section~\ref{sect Fukaya category})
of a microlocal brane $L$ does not affect calculations.

To isolate the role of the taming perturbation, let us fix a finite collection
of test branes $P_1, \cdots, P_d \subset T^*X$, and without loss of generality,
assume that they are already 
mutually transverse
and do not intersect each other or $L$ at infinity. 

\begin{lem}
Suppose the brane $L$ is equipped with two taming perturbations $\psi$ and $\psi'$.
Then for any finite collection of test objects $P_1,\ldots, P_d$ of $F(T^*X)$,
there is a quasi-isomorphism intertwining the $A_\infty$-operations
within the collection $\psi(L),  P_1,\ldots, P_d$ and the collection $\psi'(L),  P_1,\ldots, P_d$.
\end{lem}

\begin{proof}
By assumption, we have compatible almost complex structures $J$ and $J'$ such that
 the respective collections of branes 
$\psi(L),  P_0,\ldots, P_d$ 
and $\psi'(L),  P_0,\ldots, P_d$ are tame with respect to the respective induced metrics
$g$ and $g'$.

Pulling back the almost complex structures $J$ and $J'$, we obtain tame almost complex structures
$J_0=\psi^*(J)$ and $J_0'={\psi'}^*(J')$ such that
 the  branes 
$L,  P_0,\ldots, P_d$ are tame with respect to both $J_0$ and $J_0'$.
Thus we are in the setting of Lemma~~\ref{indep of ac str}, and can conclude
that the $A_\infty$-operations calculated with respect to $J_0$ and $J_0'$
are quasi-isomorphic. By construction, this is the same as a 
a quasi-isomorphism intertwining the $A_\infty$-operations
within the collection $\psi(L),  P_1,\ldots, P_d$ and the collection $\psi'(L),  P_1,\ldots, P_d$.
\end{proof}


\vskip 0.2in
\noindent
{\scriptsize
{\bf David Nadler,} Department of Mathematics, Northwestern University,
2033 Sheridan Road, Evanston, IL  60208\\
nadler@math.northwestern.edu
}


\end{document}